# Identification of Differential Equations by Dynamics-Guided Weighted Weak Form with Voting

Jiahui Cheng, Sung Ha Kang, Haomin Zhou, and Wenjing Liao *

June 4, 2025

**Abstract**

In the identification of differential equations from data, significant progresses have been made with the weak/integral formulation. In this paper, we explore the direction of finding more efficient and robust test functions adaptively given the observed data. While this is a difficult task, we propose weighting a collection of localized test functions for better identification of differential equations from a single trajectory of noisy observations on the differential equation. We find that using high dynamic regions is effective in finding the equation as well as the coefficients, and propose a dynamics indicator per differential term and weight the weak form accordingly. For stable identification against noise, we further introduce a voting strategy to identify the active features from an ensemble of recovered results by selecting the features that frequently occur in different weighting of test functions. Systematic numerical experiments are provided to demonstrate the robustness of our method.

## 1 Introduction

Recently, data-driven discovery of differential equations from experimental data has attracted much attention. Early work on automatic discovery of differential equations focused on parameter estimation [1, 2, 16, 3, 4, 19, 22], where a certain form of the underlying differential equation is known, and the goal is to identify unknown parameters using strategies such as nonlinear least squares in the differential equation.

When the form of the differential equation is unknown, the identification problem becomes more challenging since one needs to identify both the equation form and the parameters. By assuming the governing equation to be a linear combination of linear and nonlinear differential terms from a library, this identification problem can be formulated as a linear system. Among many terms, it is beneficial to find a simpler equation; hence, sparse regression is often used to choose the active features. Representative methods include Sparse Identification of Nonlinear Dynamics (SINDy) [5], Identifying Differential Equations with Numerical Time evolution (IDENT) [15], Robust IDENT [12], learning PDEs via sparse optimization [21], weak formulations in [9, 20], Weak SINDy [17, 18], WeakIdent [24], and many others [10, 11, 13, 14, 23, 25, 26].

Recently, significant progress has been made using the weak/integral formulation of differential equations [17, 18, 20, 24]. A collection of test functions are taken in the weak formulation of the PDE $\partial_t u = \sum_{l=1}^L a_l \partial_x^{\alpha_l}(u^{\beta_l})$, which gives rise to

$$\langle \partial_t u, \phi \rangle_\Omega = \sum_l a_l \langle \partial_x^{\alpha_l}(u^{\beta_l}), \phi \rangle_\Omega, \qquad (1)$$

---

*School of Mathematics, Georgia Institute of Technology, Atlanta, GA 30332-0160. Email: jcheng328@gatech.edu, kang@math.gatech.edu,hmzhou@math.gatech.edu and wliao60@gatech.edu.

where $\phi$ is a test function supported on $\Omega$ and $\langle u, v \rangle_\Omega = \int_\Omega u(x,t)v(x,t)dxdt$ is the standard inner product. By carefully selecting the test function $\phi$, the derivative acting on $u$ in (1) can be transferred to $\phi$, thereby mitigating the instability associated with numerical differentiation on noisy data. A localized test function supported on $\Omega$ enables the extraction of dynamic information from data within $\Omega$. Thus, the support and shape of the test function play a significant role in capturing relevant data features. The weak formulation proposed in [20] utilized a collection of localized test functions centered at some random points. By carefully designing the test functions, [17, 18] significantly improved the robustness of PDE identification from noisy data. The robustness was further enhanced in [24] using narrow fitting and trimming. Besides utilizing weak/integral features, Fourier features were also explored in [25], where the linear system is formed by taking the Fourier transform of the differential equation.

We explore if there is a way to design data-driven test functions to enhance the identification of differential equations. In this paper, we propose Identification of Differential Equations by Weighted weak form and Voting (IDENT-WV) to enhance the robustness and improve the accuracy of PDE identification. The proposed method utilizes a collection of dynamics-guided weighted weak forms of the differential equation, and we obtain the final result by voting. Building good weights are important and we find motivations from (i) the use of high dynamic region in WeakIdent [24]. Utilizing the high dynamic region not only helps with the efficiency of the method, but also improves the coefficient value recovery. (ii) In [25], the core regions of features are considered in the Fourier domain to utilize the meaningful data region for stable recovery. In this paper, we generalize these ideas in relation to preconditioned least-squares. Since the dynamics of the underlying PDE vary from location to location, we use different weights for localized test functions supported on different regions to highlight on high dynamic regions. To further stabilize the identification, we introduce a voting strategy based on the occurrences of active features. This voting strategy enables one to combine the identification results from multiple sets of weighted test functions to enhance robustness. Related techniques on using multiple trials of experiments can be found in Ensemble-SINDy [8]. The contributions of this paper are summarized as follows:

- We propose a weighted weak form, by combining the preconditioned least-squares with the weighted weak form for the identification of differential equations. The underlying motivation is to possibly optimize the test function for the best utilization of data in PDE identification.

- The proposed method uses occurrence voting followed by coefficient voting to stabilize the identification by removing insignificant features whose occurrence or coefficient value is low.

- We provide an error analysis for the weighted weak form to show that, by minimizing the residual in preconditioned least squares, one tends to significantly suppress the error at the support of heavily weighted test functions.

- Systematic numerical experiments are provided to demonstrate the robustness of Ident-WV in comparison with WeakIdent [24], Weak SINDy [17] and Ensemble-SINDy [8].

**Organization.** Our paper is organized as follows: We present the problem setup and the weak formulation in Section 2. Our method of Identification of Differential Equations via Weighted Weak Form with Voting (Ident-WV) is presented in Section 3. An error analysis of the weighted weak form is provided in Section 4. Our systematic numerical experiments are presented in Section 5. Finally we conclude our paper in Section 6.

## 2 Problem Setup

Let the given set of discrete and noisy observations be

$$\mathcal{D} = \{\widehat{U}_i^n = U_i^n + \epsilon_i^n\}, \text{ with } U_i^n = u(x_i, t^n) \tag{2}$$

where $\epsilon_i^n$ stands for the observation noise at $(x_i, t^n)$, which is assumed to be i.i.d. Gaussian noise. The spatial interval $[X_1, X_2]$ is partitioned into $N_x$ sub-intervals, where $x_i = X_1 + i\Delta x$, $i = 0, \ldots, N_x$, with $\Delta x = (X_2 - X_1)/N_x$ and the time interval is partitioned into $N_t$ sub-intervals such that $t^n = n\Delta t$, $n = 0, \ldots, N_t$ with $\Delta t = T/N_t$. Thus, a closed spatial-temporal domain is $[X_1, X_2] \times [0, T]$ with $X_2 > X_1$ and $T > 0$. For simplicity, we present the problem setup and our proposed algorithm when the spatial domain is one-dimensional (1D). Our method can be easily generalized to higher-dimensional spatial domains, and our numerical experiments include both 1D and 2D examples.

We assume that $U_i^n$ is a discrete sample from a governing function $u(x_i, t^n)$, and the underlying differential equation can be represented as a linear combination of linear and nonlinear terms such that

$$\partial_t u = \sum_{l=1}^{L} a_l f_l = \sum_{l=1}^{L} a_l \partial_x^{\alpha_l} T_l(u) = \sum_{l=1}^{L} a_l \partial_x^{\alpha_l} u^{\beta_l} \tag{3}$$

where $f_l = \partial_x^{\alpha_l} T_l(u) = \partial_x^{\alpha_l} u^{\beta_l}$ stands for the $l$th feature, $a_l$ is the coefficient for the $l$th feature, and $\alpha_l, \beta_l$ are nonnegative integers. We use $\partial_t, \partial_x$ to denote the partial derivative of $u$ with respect to the temporal and spatial variables respectively. Each feature $f_l$ is assumed to be the $\alpha_l$th derivative of the monomial $T_l(u) = u^{\beta_l}$, but our method is general and can be applied to other functional $T_l(u)$, such as $\sin u$, $\cos u$ and others. Let $\bar{\alpha}$ be the highest order of derivatives and $\bar{\beta}$ be the highest order of monomials in the prescribed library of features, i.e. $\bar{\alpha} = \max_l \alpha_l$, and $\bar{\beta} = \max_l \beta_l$. In this paper, we work with a general library which is defined by fixing $\bar{\alpha}$ and $\bar{\beta}$, and the prescribed feature library can be represented as $[\partial_x^{\alpha_1}(u^{\beta_1}) \quad \partial_x^{\alpha_2}(u^{\beta_2}) \quad \ldots \quad \partial_x^{\alpha_L}(u^{\beta_L})] \in \mathbb{R}^{1 \times L}$. The PDE in (3) is expressed as

$$\partial_t u = [\partial_x^{\alpha_1}(u^{\beta_1}) \quad \partial_x^{\alpha_2}(u^{\beta_2}) \quad \ldots \quad \partial_x^{\alpha_L}(u^{\beta_L})] \boldsymbol{a}, \text{ where } \boldsymbol{a} = (a_1, \ldots, a_L)^\top \in \mathbb{R}^L \tag{4}$$

represents the unknown coefficients in the governing function. While the library of features may have a large number of terms, the coefficient vector $\boldsymbol{a}$ in (4) is usually assumed to be sparse, with a small number of nonzero entries, since we are interested in simple equation representing physical world phenomenon. The features associated with the nonzero coefficients are called active features in the governing function, driving the dynamics of the PDE.

To utilize the weak form, we define the test function in the following way. Let $\Omega \subset [X_1, X_2] \times [0, T]$ be a spatial-temporal sub-domain of $[X_1, X_2] \times [0, T]$, and $\phi : \Omega \mapsto \mathbb{R}$ be a test function supported on $\Omega$. The test function $\phi$ and its spatial derivatives up to order of $\bar{\alpha}$ vanish on boundary of $\Omega$. Taking an inner product between the test function $\phi$ and (3) gives rise to

$$\langle \partial_t u, \phi \rangle_\Omega = \sum_l a_l \langle \partial_x^{\alpha_l}(u^{\beta_l}), \phi \rangle_\Omega,$$

as in (1), and by integration by parts of (1) gives rise to

$$-\langle u, \partial_t \phi \rangle_\Omega = \sum_l a_l (-1)^{|\alpha_l|} \langle u^{\beta_l}, \partial_x^{\alpha_l} \phi \rangle_\Omega. \tag{5}$$

The weak formulation transfers the derivatives from $u$ to the test function $\phi$, which gives rises to a better stability in handling high order derivatives for noisy data. We use a collection of test basis functions $\{\phi_h\}_{h=1}^H$, where each $\phi_h$ is localized and supported on a local region denoted by $\Omega_h$ as in [17, 24]. The weak formulation in (5) in discrete form yields the linear system:

$$W\boldsymbol{a} = \boldsymbol{b} \tag{6}$$

where $W \in \mathbb{R}^{H \times L}$ and $\boldsymbol{b} \in \mathbb{R}^H$ are matrices and vectors containing the features in the weak form:

$$W = \begin{pmatrix} | & | & \cdots & | \\ (-1)^{|\alpha_1|}\langle u^{\beta_1}, \partial_x^{\alpha_1}\phi_h\rangle_{\Omega_h} & (-1)^{|\alpha_2|}\langle u^{\beta_2}, \partial_x^{\alpha_2}\phi_h\rangle_{\Omega_h} & \cdots & (-1)^{|\alpha_L|}\langle u^{\beta_L}, \partial_x^{\alpha_L}\phi_h\rangle_{\Omega_h} \\ | & | & \cdots & | \end{pmatrix} \in \mathbb{R}^{H \times L},$$

$$\boldsymbol{b} = \begin{pmatrix} | \\ -\langle u, \partial_t\phi_h\rangle_{\Omega_h} \\ | \end{pmatrix} \in \mathbb{R}^H, \tag{7}$$

We refer to $W = (\boldsymbol{w}_1, \boldsymbol{w}_2, \ldots, \boldsymbol{w}_L)$ as the feature matrix representing the feature library, where $\boldsymbol{w}_l$ denotes the $l$th column of $W$.

The objective of this paper is to identify the underlying differential equation by identifying both the support and coefficient values of $\boldsymbol{a}$ in (6) from the single observation of a given discrete and noisy data set $\mathcal{D}$ as in (2).

For notation, we use bold letters to denote vectors and regular letters to denote matrices and scalars. For a row (or column) vector $\boldsymbol{a} = (a_1, \ldots, a_L) \in \mathbb{R}^L$, we denote its support by $\text{supp}(\boldsymbol{a}) = \{l : a_l \neq 0\}$. We denote $\text{diag}(\boldsymbol{a}) \in \mathbb{R}^{L \times L}$ as the diagonal matrix whose diagonal entries are $\boldsymbol{a}$. We use $\#A$ to denote the cardinality of the set $A$. If $\mathcal{B}$ is an index set and $W$ is a matrix, $W_\mathcal{B}$ denotes the submatrix of $W$ where the column index is restricted to $\mathcal{B}$. We use $\mathbb{1}$ to denote the indicator function.

## 3 Identification of Differential Equations via Dynamics-guided Weighted Weak Form with Voting (Ident-WV)

Since using test functions and formulating differential equations in a weak form demonstrate various advantages, we take this idea further. (1) We first introduce the use of dynamically guided weighted test functions to improve identification accuracy. We find that using high dynamic regions is effective for identifying both the differential equation and its coefficients, and we propose a dynamics indicator for each differential term and apply it to weight the weak form. (2) To further improve the stability against noisy data, we incorporate occurrence voting followed by coefficient voting for equation identification. In this section, we present the details of both the dynamics-guided weighted weak form and the occurrence voting, then we illustrate this procedure carefully with a numerical example, before we present the error analysis in the following section.

### 3.1 Dynamics-guided weighted weak form

Our main idea is to learn/find the weights of test functions in the weak form for more robust identification. In WeakIdent [24], the high dynamic region of the PDE defined according to the reference feature $(u^2)_x$ is used to enhance the coefficient recovery. We generalize the idea of exploiting high dynamic regions defined according to more reference features and include the time derivative.

For each reference feature $g_m = \partial_t^\gamma \partial_x^\alpha u^\beta$, we define the **dynamics indicator** associated with the feature $\partial_t^\gamma \partial_x^\alpha u^\beta$ as the leading error coefficient in the noise expansion. A Taylor expansion with respect to $\varepsilon$ gives rise to

$$\int_{\Omega_h} \partial_t^\gamma \partial_x^\alpha (u+\varepsilon)^\beta \phi_h dxdt = \int_{\Omega_h} \partial_t^\gamma \partial_x^\alpha \left( u^\beta + \binom{\beta}{1} u^{\beta-1}\varepsilon + O(\varepsilon^2) \right) \phi_h dxdt$$

$$= \int_{\Omega_h} \partial_t^\gamma \partial_x^\alpha (u^\beta) \phi_h dxdt + \left( \beta \int_{\Omega_h} \partial_t^\gamma \partial_x^\alpha (u^{\beta-1}) \phi_h dxdt \right) \varepsilon + O(\varepsilon^2).$$

The leading error coefficient associated with the reference feature $\partial_t^\gamma \partial_x^\alpha u^\beta$ and the test function $\phi_h$ becomes

$$r(h, \partial_t^\gamma \partial_x^\alpha u^\beta) := \beta \left| \int_{\Omega_h} \partial_t^\gamma \partial_x^\alpha (u^{\beta-1}) \phi_h dxdt \right| = \beta \left| \int_{\Omega_h} u^{\beta-1} \partial_t^\gamma \partial_x^\alpha (\phi_h) dxdt \right|, \ h=1,2,...,H, \quad (8)$$

where the second equality follows from our choice of test functions such that $\phi_h$ and its spatial derivatives up to order $\alpha - 1$ and time derivatives up to order $\gamma - 1$ all vanish on the boundary of $\Omega_h$. A large dynamics indicator occurs when the reference feature has a large leading error coefficient, and vice versa. As shown in Figure 2 in Section 3.4, the high/low dynamic regions of the PDE can be partially captured by the dynamic indicator defined in (8). When $\beta = 1$, the dynamic indicator in (8) becomes independent of $u$, since $\left| \int_{\Omega_h} \partial_t^\gamma \partial_x^\alpha (\phi_h) dxdt \right|$ is the same for all $h$'s. In other words, when $\beta = 1$, the dynamic indicator associated with $u$ is constant for all $h$'s, giving rise to equally weighted test functions.

For a fixed reference feature $g_m = \partial_t^\gamma \partial_x^\alpha u^\beta$, the dynamic indicator $r(h, g_m) = r(h, \partial_t^\gamma \partial_x^\alpha u^\beta)$ depends on the row index $h$ (test function index) and the reference feature $g_m$. For this fixed reference feature, we put the dynamics indicator for each row together as a column vector and define **a weight matrix** $R^{(m)}$:

$$R^{(m)} = R[g_m] = R[\partial_t^\gamma \partial_x^\alpha u^\beta] = \mathrm{diag}(\boldsymbol{r}(\partial_t^\gamma \partial_x^\alpha u^\beta)) \text{ where } \boldsymbol{r}(\partial_t^\gamma \partial_x^\alpha u^\beta) = \begin{pmatrix} | \\ r(h, \partial_t^\gamma \partial_x^\alpha u^\beta) \\ | \end{pmatrix} \in \mathbb{R}^H. \quad (9)$$

We define the **Dynamics-Guided Weighted Weak Form** as:

$$R^{(m)} W \boldsymbol{a} = R^{(m)} \boldsymbol{b}, \quad (10)$$

where $W \in \mathbb{R}^{H \times L}$ and $\boldsymbol{b} \in \mathbb{R}^H$ are defined in (7). In (10), we use the diagonal matrix $R^{(m)} = R[g_m] = R(\partial_t^\gamma \partial_x^\alpha u^\beta)$ to emphasize the high dynamic regions defined according to the reference feature $g_m$. Next, we use WeakIdent [24] to recover the sparse coefficient $\boldsymbol{a}$ from the linear system in (10). We denote

$$\widehat{\boldsymbol{a}}^{(m)} = \mathrm{WeakIdent}(R^{(m)} W, R^{(m)} b) \quad (11)$$

as the recovered coefficient vector by the dynamics-guided weak form with weight matrix $R^{(m)}$. The recovered support for active features is given by $\mathrm{supp}(\widehat{\boldsymbol{a}}^{(m)})$.

In this paper, for PDEs in the form of (3), we use the weight matrices $R$ from reference features which contain the spatial derivatives of the dependent variable (e.g. $u^2$) up to second order, and its first-order time derivative, i.e. $u$, $u_x$, $u_{xx}$, $u_t$, $u^2$, $(u^2)_x$, $(u^2)_{xx}$, $(u^2)_t$. These eight features are given by $u$ or $u^2$, or their spatial derivatives up to the second order, or their first-order time derivative. According to (8), the dynamics indicators of $u$, $u_x$, $u_{xx}$, $u_t$ are the same, so only pick $u$ as a representative reference feature for $u$, $u_x$, $u_{xx}$, $u_t$. Hence, eight reference features are reduced to five reference features $g_1 = u$, $g_2 = u^2$, $g_3 = (u^2)_x$, $g_4 = (u^2)_{xx}$, $g_5 = (u^2)_t$. These five reference features give rise to five weight matrices $R^{(m)} = R[g_m]$ for $m = 1, \ldots, M$. This includes the high dynamic region of $(u^2)_x$ used in WeakIdent [24], and also includes $u$, $u^2$, $(u^2)_{xx}$ and $(u^2)_t$.

## 3.2 Occurrence and coefficient voting

The dynamic-guided weighted weak form in (10) for each reference feature $g_m, m = 1, \ldots, M$ gives rise to a coefficient recovery for each $m$. We also add occurrence and coefficient voting to the collection of results given by (10) for $m = 1, \ldots, M$.

By using the dynamics indicators, we have a collection of $M$ weight matrices $R^{(m)}$, $m = 1, \ldots, M$, highlighting high dynamic regions of the PDE. Suppose that the true coefficient vector $\boldsymbol{a}$ has sparsity $S$ such that the underlying PDE has $S$ active features. Then, sparse regression

$$\boldsymbol{a}^{(m)} = \operatorname*{argmin}_{\boldsymbol{z} \in \mathbb{R}^L : \|\boldsymbol{z}\|_0 \leq S} \|R^{(m)}W\boldsymbol{z} - R^{(m)}\boldsymbol{b}\|_2$$

with the preconditioner $R^{(m)}$ allows one to suppress the residual more in the rows where $R^{(m)}$ has larger magnitudes. This yields a good estimate of the coefficient vector if $R^{(m)}W$ satisfies the incoherence property [7] or the restricted isometry property [6]. We propose the **occurrence voting** of each feature, given by

$$\text{occurrence}_l = \frac{1}{M} \sum_{m=1}^{M} \mathbb{1}\left\{ a_l^{(m)} \neq 0 : \boldsymbol{a}^{(m)} = \operatorname*{argmin}_{\boldsymbol{z} \in \mathbb{R}^L : \|\boldsymbol{z}\|_0 \leq S} \|R^{(m)}W\boldsymbol{z} - R^{(m)}\boldsymbol{b}\| \right\} \quad (12)$$

to determine active features. In practice, we do not have prior information about the true sparsity $S$ and the sparse regression problem in (12) is combinatorially hard to solve. Instead, we approximately find each sparse solution $\widehat{\boldsymbol{a}}^{(m)}$ by WeakIdent [24] and vote for the active features based on the occurrence of each feature.

In addition, we propose the **coefficient voting**. First, active features are refined so that insignificant features with small coefficient values are removed. The average coefficient value for the active features in $\mathcal{B}$ is computed as,

$$\bar{a}_l = \frac{1}{M} \sum_{m=1}^{M} |\widehat{a}_l^{(m)}|, \quad \text{for } l \in \mathcal{B}.$$

We remove the features with relatively small coefficients from the support $\mathcal{B}$ to update the support to:

$$\mathcal{C} = \{l \in \mathcal{B} : \frac{\bar{a}_l}{\max_l \bar{a}_l} \geq \upsilon\}, \quad (13)$$

where $\upsilon$ is a preset coefficient threshold, which is set to be 5% in this paper.

## 3.3 Coefficient recovery

After occurrence and coefficient voting gives rise to the recovered support $\mathcal{C}$ in (13), we solve the linear system (6) rescaled by the average leading error coefficient [17, 24]:

$$e\langle l \rangle = \frac{1}{H} \sum_{h=1}^{H} e(h,l), \text{ where } e(h,l) := \beta_l \left| \int_{\Omega_h} u^{\beta_l - 1} \partial_x^{\alpha_l} \phi_h dx dt \right|, \quad (14)$$

for $l = 1, \ldots, L$. Here $e(h,l)$ represents the leading error coefficient for the $h$-th test function $\phi_h$ and the $l$-th feature $f_l$. The $e\langle l \rangle$ in (14) averages over the leading error coefficient for all test function indexed by $h = 1, \ldots, H$.

To obtain the coefficient vector $\widehat{\boldsymbol{a}}$ supported on $\mathcal{C}$, we rescale the weak feature matrix $W$ to $\widetilde{W}$ such that

$$\widetilde{W}\widetilde{\boldsymbol{a}} = \boldsymbol{b} \quad (15)$$

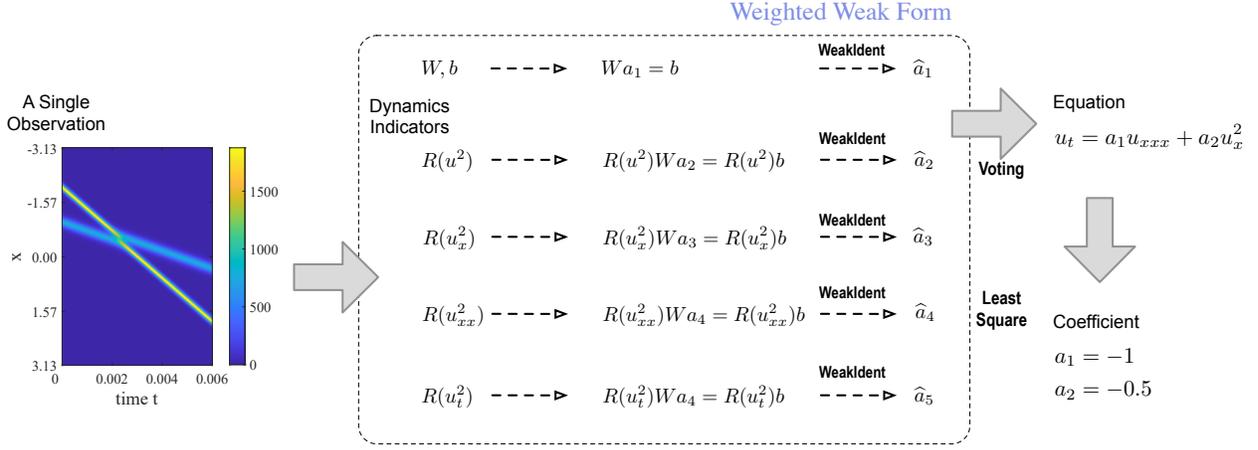

Figure 1: The workflow of Ident-WV applied to the KdV equation $u_t = -u_{xxx} - 0.5(u^2)_x$.

where $\widetilde{D} = \mathrm{diag}(e\langle 1\rangle, , e\langle 2\rangle, \ldots, e\langle L\rangle) \in \mathbb{R}^{L\times L}$, $\widetilde{\boldsymbol{a}} = \widetilde{D}\boldsymbol{a} \in \mathbb{R}^L$, and $\widetilde{W} = W\widetilde{D}^{-1}$. The recovered coefficient $\widehat{\boldsymbol{a}}$ is computed through least squares:

$$\widehat{\boldsymbol{a}} = \widetilde{D}^{-1}\widetilde{\boldsymbol{a}} \in \mathbb{R}^L, \quad \text{where} \quad \widetilde{\boldsymbol{a}} = \underset{\mathrm{Supp}(\boldsymbol{z})\subset \mathcal{C}}{\mathrm{argmin}} \left\|\widetilde{W}\boldsymbol{z} - \boldsymbol{b}\right\|. \tag{16}$$

By rescaling the linear system by the average leading error coefficient, one can significantly improve the condition number of the least squares problem and therefore enhance the stability of coefficient recovery [24]. Compared to [24], IDENT-WV does not utilize the narrow fit technique in [24] where the rows in the high dynamic region of $(u^2)_x$ are selected for coefficient recovery.

### 3.4 Numerical illustration of the proposed method

In this section, we summarize the process of the proposed Ident-WV. As an example, we consider a Korteweg–De Vries equation (KdV) equation:

$$u_t = -u_{xxx} - 0.5(u^2)_x \tag{17}$$

from a single observation of the solution with 30% noise as in (26). The initial condition is $u(x,0) = 3 \times 25^2 \times \mathrm{sech}(0.5 \times 25(x+2.0))^2 + 3 \times 16^2 \times \mathrm{sech}(0.5 \times 16(x+1.0))^2$. Figure 1 demonstrates the workflow of our proposed Ident-WV algorithm.

First, we construct $W, b$ and five dynamics indicators $R^{(1)} = R[u]$, $R^{(2)} = R[u^2]$, $R^{(3)} = R[(u^2)_x]$, $R^{(4)} = R[(u^2)_{xx}]$, $R^{(5)} = R[(u^2)_t]$ from noisy data. Even with a large amount of noise on the given data, the dynamics indicators produce good representations of the high/low dynamic regions of the solution, which are shown in Figure 2. The dynamics indicators associated with different reference features give rise to different emphases on the high dynamic regions, as represented as yellow regions in Figure 2. For example, if we take $u^2$ as the reference feature, the dynamics indicator is large when leading error coefficient of $u^2$ has a large magnitude.

Second, we perform dynamic-guided weighted weak form to predict the differential equation for each dynamic indicator and collect the occurrence and average magnitude of each feature. Figure 3 illustrates the occurrence of recovered features from dynamic-guided weighted weak form with five reference features, where we apply an occurrence voting followed by a coefficient voting to select features. The occurrences of true features are 3/5 and 4/5, and some false features are identified with the occurrence 1/5 or 2/5. Our first vote is based on the feature occurrence, and our second

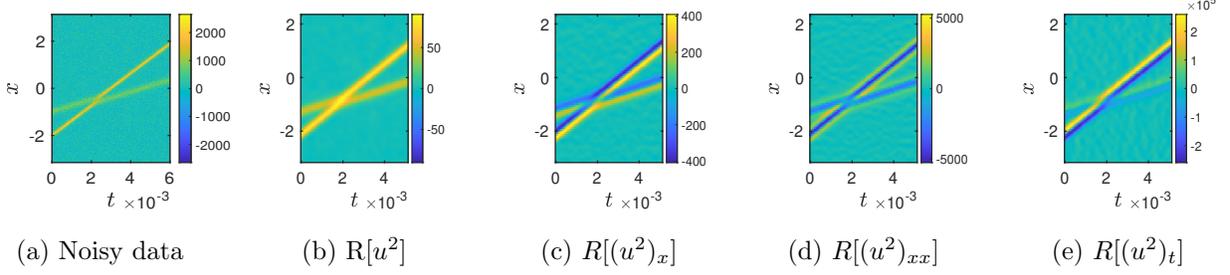

(a) Noisy data  (b) R[$u^2$]  (c) $R[(u^2)_x]$  (d) $R[(u^2)_{xx}]$  (e) $R[(u^2)_t]$

Figure 2: [First Step of Ident-WV] Dynamics indicators for the KdV equations. The first figure (a) displays the noisy data. (b) - (e) display the dynamics indicators $R[u^2], R[(u^2)_x], R[(u^2)_{xx}], R[(u^2)_t]$ when the reference feature is $u^2, u_x^2, u_{xx}^2, u_t^2$ respectively. In (b) - (e), the pixel value at $(x^{(h)}, t^{(h)})$ stands for the dynamics indicator $r(h, \partial_t^\gamma \partial_x^\alpha u^\beta)$ defined in (8) normalized by the maximum magnitude of $r(h, \partial_t^\gamma \partial_x^\alpha u^\beta)$.

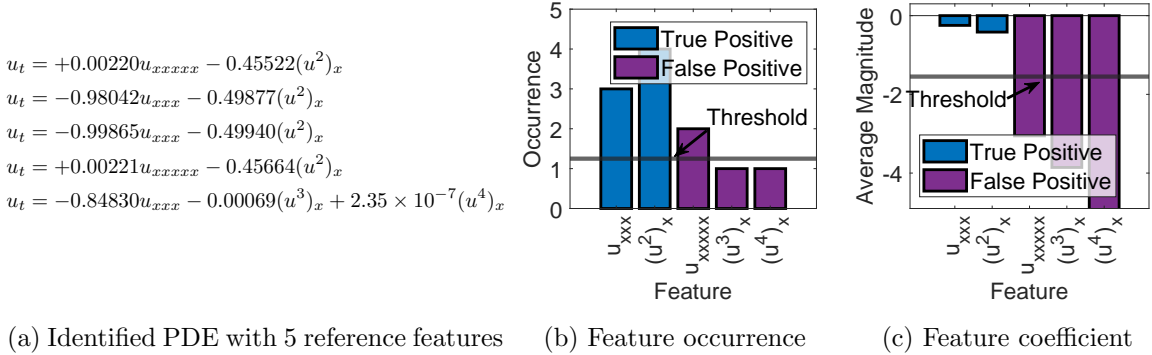

$u_t = +0.00220 u_{xxxxx} - 0.45522(u^2)_x$
$u_t = -0.98042 u_{xxx} - 0.49877(u^2)_x$
$u_t = -0.99865 u_{xxx} - 0.49940(u^2)_x$
$u_t = +0.00221 u_{xxxxx} - 0.45664(u^2)_x$
$u_t = -0.84830 u_{xxx} - 0.00069(u^3)_x + 2.35 \times 10^{-7}(u^4)_x$

(a) Identified PDE with 5 reference features   (b) Feature occurrence   (c) Feature coefficient

Figure 3: [Second Step of Ident-WV] (a) shows the identified differential equations by the weighted weak forms with 5 reference features: $u, u^2, (u^2)_x, (u^2)_{xx}, (u^2)_t$. (b) Occurrence voting and (c) Coefficient Voting. The true features are marked in blue and false features are in purple.

vote is based on the coefficient value. Let $\rho$ and $\nu$ be two thresholding parameters. We remove the features whose occurrence is low (below $\rho$), after which we trim the features whose average amplitude is low (below $\nu$ of the maximum amplitude). The threshold $\rho = 20\%$ and $\nu = 5\%$ is represented by the horizontal line in Figure 3. In this example, Ident-WV identifies the correct features of the KdV equation. Finally, the coefficient vector is computed by a least squares fitting of the rescaled linear system in (16). The algorithm is summarized in Algorithm 1.

## 4  An Error Analysis of Weighted Weak Form

Since each test function $\phi_h$ is locally supported at $\Omega_h$, introducing the weight parameters allows us to emphasize more on the local region with a larger weight. We present an error analysis to show that the least square residual is a weighted sum of the error on the $\Omega_h$'s. Thus, minimizing the residual suppresses the error more heavily on the high dynamic region where the test function has a larger weight.

Let the true coefficient vector for the underlying PDE (5) be $\boldsymbol{a}^*$ such that the noiseless data $\{u(x_i, t^n)\}$ simulated from the PDE with the coefficient vector $\boldsymbol{a}^*$. Denote the true support as $\text{supp}^* = \text{supp}(\boldsymbol{a}^*)$. When we formulate the discrete linear system $RW\boldsymbol{a} = R\boldsymbol{b}$ in (10) from the

**Algorithm 1** Ident-WV

**Require:** $W \in \mathbb{R}^{H \times L}, \boldsymbol{b} \in \mathbb{R}^H$ in (10); A list of $M$ reference features with dynamics indicators $\boldsymbol{r}_m : m = 1, \ldots, M\}$; Parameters $\rho = 25\%$ and $\upsilon = 5\%$.
**Ensure:** Recovered coefficient $\widehat{\boldsymbol{a}}$.
1: Initialize the occurrence of each feature: $n_l = 0$ for $l = 1, \ldots, L$.
2: **for** $m = 1, \ldots, M$ **do**
3:     $R^{(m)} = \text{diag}(\boldsymbol{r}_m)$.     ▷ Construct the row multiplication matrix.
4:     $\widehat{\boldsymbol{a}}^{(m)} = \text{WeakIdent}(W, R^{(m)}, \boldsymbol{b})$.
5:     **for** $l = 1, \ldots, L$ **do**
6:         **if** $\widehat{a}_l^{(m)} \neq 0$ **then** $n_l = n_l + 1$.     ▷ Accumulate the occurrence.
7: Occurrence$_l = n_l/M$ for $l = 1, \ldots, L$.
8: **First Voting**: Trim features with small occurrence to the support $\mathcal{B} = \{l : \text{occurrence}_l \geq \rho\}$.
9: Average coefficient: $\bar{a}_l = \frac{1}{M}(\sum_{m=1}^M |\widehat{a}_l^{(m)}|)$ for $l \in \mathcal{B}$, and $\bar{a}_l = 0$, for $l \in \mathcal{B}^\complement$.
10: **Second Voting**: Trim features with small coefficients to the support $\mathcal{C} = \{l \in \mathcal{B} : \frac{\bar{a}_l}{\max_l \bar{a}_l} \geq \upsilon\}$.
11: **Recovered Coefficient**: $\widehat{\boldsymbol{a}}$ given by (16).

noisy data in $\mathcal{D}$ in (2), we can express the residual error as

$$\boldsymbol{e} = \boldsymbol{e}_{\boldsymbol{a}^*} + \boldsymbol{e}_{\boldsymbol{a}}, \text{ with } \boldsymbol{e}_{\boldsymbol{a}^*} := RW\boldsymbol{a}^* - R\boldsymbol{b} \text{ and } \boldsymbol{e}_{\boldsymbol{a}} := RW(\boldsymbol{a} - \boldsymbol{a}^*). \quad (18)$$

where the error $\boldsymbol{e}_{\boldsymbol{a}^*}$ arises from numerical integration as well as noise even when the true coefficient vector $\boldsymbol{a}^*$ is put in the weighted weak form, and the error $\boldsymbol{e}_{\boldsymbol{a}}$ quantifies the coefficient matching error.

We first analyze the $\boldsymbol{e}_{\boldsymbol{a}^*}$ error. For the $h$-th weighted test function $r_h \phi_h$ supported on the domain $\Omega_h$, the $h$-th entry in $\boldsymbol{e}_{\boldsymbol{a}^*}$ error

$$(\boldsymbol{e}_{\boldsymbol{a}^*})_h = r_h \sum_{(x_j, t^k) \in \Omega_h} \widehat{U}_j^k \partial_t \phi_h(x_j, t^k) \Delta x \Delta t + r_h \sum_{l \in \text{supp}^*} a_l^* (-1)^{|\alpha_l|} \sum_{(x_j, t^k) \in \Omega_h} (\widehat{U}_j^k)^{\beta_l} \partial_x^{\alpha_l} \phi_h(x_j, t^k) \Delta x \Delta t. \quad (19)$$

By defining a point-wise residual as

$$\mathcal{R}(u, x, t) = u(x, t) \partial_t \phi(x, t) + \sum_{l \in \text{supp}^*} a_l^* (-1)^{|\alpha_l|} [u(x, t)]^{\beta_l} \partial_x^{\alpha_l} \phi(x, t), \quad (20)$$

the error $(\boldsymbol{e}_{\boldsymbol{a}^*})_h$ in (19) can be written as

$$(\boldsymbol{e}_{\boldsymbol{a}^*})_h = \sum_{(x_j, t^k) \in \Omega_h} r_h \mathcal{R}(\widehat{U}, x_j, t^k) \Delta x \Delta t. \quad (21)$$

The following theorem provides an upper bound on the $\boldsymbol{e}_{\boldsymbol{a}^*}$ error, to quantify the effects from numerical integration and noise, even when the true coefficient vector $\boldsymbol{a}^*$ is put in the weighted weak form.

**Theorem 1.** *Consider the differential equation in (3) with the true coefficient vector $\boldsymbol{a}^*$. Assume that the noise $\epsilon_i^k$ in (2) are zero-mean i.i.d. bounded random variables such that $|\epsilon_i^k| \leq \epsilon$ for some $0 < \epsilon \ll 1$. The residual error of the discrete linear system $\boldsymbol{e}_{\boldsymbol{a}^*}$ for the true coefficient vector $\boldsymbol{a}^*$ satisfies*

$$|(\boldsymbol{e}_{\boldsymbol{a}^*})_h| \leq O[(\Delta x \Delta t)^2] + r_h S_h^* \epsilon + \mathcal{O}(r_h \epsilon^2), \ h = 1, \ldots, H$$

*with $S_h^* = \sum_{(x_j, t^k) \in \Omega_h} \left| \partial_t \phi(x_j, t^k) + \sum_{l \in \text{supp}^*} a_l^* (-1)^{|\alpha_l|} \beta_l (U_j^k)^{\beta_l - 1} \partial_x^{\alpha_l} \phi(x_j, t^k) \right| \Delta x \Delta t.*

*Proof.* The error $e_{a^*}$ can be decomposed into two terms representing the error arising from noise and numerical integration, respectively:

$$e_{a^*} = e_{a^*}^{\text{noise}} + e_{a^*}^{\text{int}}$$

where

$$(e_{a^*}^{\text{noise}})_h = r_h \sum_{(x_j,t^k)\in\Omega_h} \mathcal{R}(\widehat{U}, x_j, t^k)\Delta x \Delta t - r_h \sum_{(x_j,t^k)\in\Omega_h} \mathcal{R}(U, x_j, t^k)\Delta x \Delta t,$$

$$(e_{a^*}^{\text{int}})_h = r_h \sum_{(x_j,t^k)\in\Omega_h} \mathcal{R}(U, x_j, t^k)\Delta x \Delta t - r_h \int_{\Omega_h} \mathcal{R}(u, x, t) dx dt.$$

In this decomposition, $\widehat{U}$ and $U$ represent noisy and noiseless data respectively. The $e_{a^*}^{\text{noise}}$ term represents the error from noise. The $h$th entry of $e_{a^*}^{\text{noise}}$ can be expressed as

$$\begin{aligned}
(e_{a^*}^{\text{noise}})_h &= r_h \sum_{(x_j,t^k)\in\Omega_h} \mathcal{R}(\widehat{U}, x_j, t^k)\Delta x \Delta t - r_h \sum_{(x_j,t^k)\in\Omega_h} \mathcal{R}(U, x_j, t^k)\Delta x \Delta t, \\
&= r_h \Delta x \Delta t \sum_{(x_j,t^k)\in\Omega_h} \left( \mathcal{R}(\widehat{U}, x_j, t^k) - \mathcal{R}(U, x_j, t^k) \right) \\
&= r_h \Delta x \Delta t \sum_{(x_j,t^k)\in\Omega_h} \left( \epsilon_j^k \partial_t \phi(x_j, t^k) + \sum_{l\in\text{supp}^*} a_l^*(-1)^{|\alpha_l|}\epsilon_j^k \left( \sum_{r=1}^{\beta_l} \binom{\beta_l}{r} (\epsilon_j^k)^{r-1}(U_j^k)^{\beta_l - r} \right) \partial_x^{\alpha_l}\phi(x_j, t^k) \right) \\
&= r_h \Delta x \Delta t \sum_{(x_j,t^k)\in\Omega_h} \left( \partial_t\phi(x_j, t^k) + \sum_{l\in\text{supp}^*} a_l^*(-1)^{|\alpha_l|}\beta_l(U_j^k)^{\beta_l - 1}\partial_x^{\alpha_l}\phi(x_j, t^k) \right)\epsilon_j^k + \mathcal{O}((\epsilon_j^k)^2)
\end{aligned}$$

Thus, we obtain

$$|(e_{a^*}^{\text{noise}})_h| \leq r_h S_h^* \epsilon + \mathcal{O}(r_h \epsilon^2). \tag{22}$$

For the numerical integration error, the $h$th entry is

$$(e_{a^*}^{\text{int}})_h = r_h \sum_{(x_j,t^k)\in\Omega_h} \mathcal{R}(U, x_j, t^k)\Delta x \Delta t - r_h \int_{\Omega_h} \mathcal{R}(u, x, t) dx dt = r_h \sum_{(x_j,t^k)\in\Omega_h} \mathcal{R}(U, x_j, t^k)\Delta x \Delta t, \tag{23}$$

the second equality in (23) holds since the true equation satisfies $\int_{\Omega_h} \mathcal{R}(u, x, t) dx dt = 0$. The numerical integration is carried by the trapezoidal rule, which gives rise to the second order error:

$$|(e_{a^*}^{\text{int}})_h| \leq C(\Delta x \Delta t)^2 \max\left( \sup_{(x,t)\in\Omega_h} \left|\frac{\partial^2 \mathcal{R}}{\partial x^2}(u,x,t)\right|, \sup_{(x,t)\in\Omega_h} \left|\frac{\partial^2 \mathcal{R}}{\partial t^2}(u,x,t)\right|, \sup_{(x,t)\in\Omega_h} \left|\frac{\partial^2 \mathcal{R}}{\partial x \partial t}(u,x,t)\right| \right), \tag{24}$$

where $C$ is a constant depending on the volume of $\Omega_h$. Combining (22) and (24) gives the result. □

Theorem 1 quantifies the effects from noise and numerical integration of the dynamics-guided weighted weak form. First of all, it shows that the discrete linear system in (10) is consistent for the true coefficient vector $a^*$ such that

$$\lim_{\Delta x, \Delta t, \epsilon \to 0} e_{a^*} = \mathbf{0}.$$

Secondly, Theorem 1 shows that the weighted weak form is robust to noise. The weighted weak form enjoys the same robustness as the weak/integral form, which is a special case with $r_h = 1$ for all $h$. The numerical integration error scales quadratically with the grid spacing $\Delta x \Delta t$ and the error from noise scales linearly with the noise level $\epsilon$. When data are noisy, the error in the weighted weak form is significantly smaller than that in the differential form. In comparison, it was shown in [15] that the error for the discretized system under the differential form is on the order of

$$O(\Delta t + \Delta x^{p+1-r} + \frac{\epsilon}{\Delta t} + \frac{\epsilon}{\Delta x^r}) \qquad (25)$$

where $r$ is the highest order of derivatives for the features in the true support, and the numerical differentiation is carried by interpolating the data by a $p$th order polynomial. As $\Delta x$ and $\Delta t$ decrease to 0, the error in (25) for the differential form may blow up, while the error for the dynamics-guided weighted weak form converges to zero.

In addition to $e_{a^*}$, the residual error for any coefficient vector $a$ in (18) has another error term $e_a = RW(a - a^*)$, which measures the mismatch error between $a$ and $a^*$. In the dynamics-guided weighted weak form, the residual is a weighted sum of the local errors on the $\Omega_h$'s. The weight matrix serves like a pre-conditioner in solving this pre-conditioned least squares problem, which tends to suppress the residual at the support of heavily weighted test functions.

## 5 Numerical Experiments

In this section, we present various numerical results to demonstrate the robustness of Ident-WV. To quantify the noise level of data, we use the Noise-to-Signal Ratio (NSR), denoted by $\sigma_{\text{NSR}}$ such that

$$\sigma^2_{\text{NSR}} = \frac{\sigma^2}{\frac{1}{\mathbb{N}_t \mathbb{N}_x} \sum_{i,n} |U_i^n - (\max_{i,n} U_i^n + \min_{i,n} U_i^n)/2|^2} \qquad (26)$$

for i.i.d. Gaussian noise such that $\epsilon_i^n \sim \mathbb{N}(0, \sigma^2)$ as in [24]. This NSR definition follows that in WeakIdent [24], and differs from the definitions in Ident [15] and WeakSINDy [17]. This definition has the advantage of being invariant to the shift of $U$ such that the data sets $\{U_i^n + c\}$ and $\{U_i^n\}$ have the same NSR when $\sigma$ is fixed.

To compare the results, we use the following three measures. The first two, the True Positive Ratio (TPR) and the Positive Predictive Value (PPV), measure the accuracy of the coefficient support identification that these measure the accuracy of the form of the equations. Denote $a \in \mathbb{R}^L$ and $\widehat{a} \in \mathbb{R}^L$ as the true coefficient vector and the recovered coefficient vector, respectively. The True Positive Ratio (TPR) is defined as the ratio between the cardinality of the correctly identified features and the cardinality of the true support,

$$\text{TPR} = \frac{\#\{l : a_l \neq 0 \text{ and } \widehat{a}_l \neq 0\}}{\#\{l : a_l \neq 0\}}. \qquad (27)$$

The TPR is equal to 1 if all the true features are identified. When the TPR is below 1, only a fraction of the true features are identified, and some true features are missing in the recovery. The Positive Predictive Value (PPV) is defined as the ratio between the cardinality of the correctly identified support and the cardinality of the recovered support,

$$\text{PPV} = \frac{\#\{l : a_l \neq 0 \text{ and } \widehat{a}_l \neq 0\}}{\#\{l : \widehat{a}_l \neq 0\}}. \qquad (28)$$

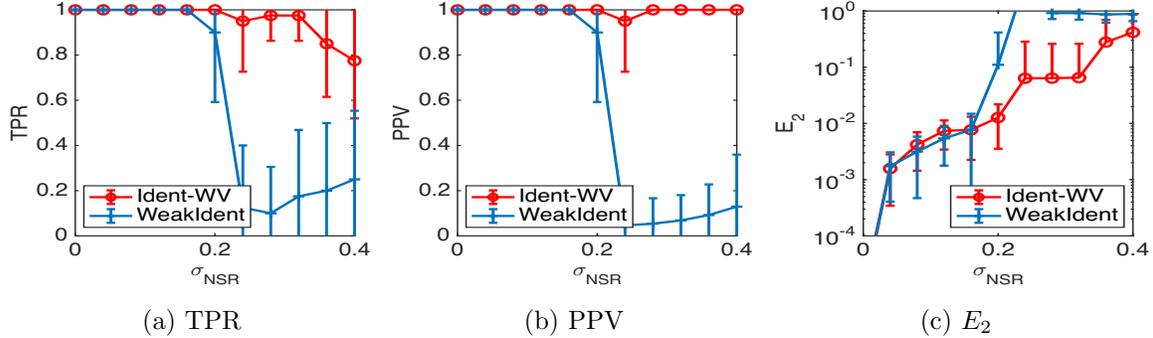

(a) TPR  (b) PPV  (c) $E_2$

Figure 4: KdV equation (17) shown in Section 3.4. We show the comparison against WeakIdent [24], via the average TPR (27), PPV (28) and $E_2$ error. The curve represents the average error and the standard deviation is represented by the vertical bar. Ident-WV (red) shows improved robustness compared to WeakIdent (blue), when the noise level is higher than 20%, showing higher TPR and PPV, and smaller $E_2$ coefficient error.

The PPV is equal to 1 if the recovered support is a subset of the true support. When the PPV is below 1, it indicates the presence of false positive features. The third measure is the relative $\ell_2$ coefficient error:

$$E_2 = \|\boldsymbol{a} - \widehat{\boldsymbol{a}}\|_2 / \|\boldsymbol{a}\|_2, \tag{29}$$

which measures the accuracy of the coefficient recovery.

In the following subsections, we present systematic numerical experiments, comparing the proposed IDENT-WV with WeakIdent [24], WeakSINDy [17] and Ensemble-SINDy [8]. We first present the identification result for the KdV equation (17) in Subsection 5.1, and then some 1D equations with initial condition at various frequencies in Subsection 5.2. The results for more 1D equations with higher order derivatives and 2D equations are given in Subsection 5.3. Finally, we compare Ident-WV and Ensemble-SINDy [8] in Subsection 5.4. For a fixed NSR, we run 20 experiments and average the error.

## 5.1 First example: the KdV equation (17) in comparison with WeakIdent

In Section 3.4, we have demonstrated the identification procedure for the KdV equation $u_t = -u_{xxx} - 0.5u_x^2$ with 30% noise. Figure 4 further shows the statistical results for 20 experiments per each noise level varying from 0% to 40%. Figures 4 (a) and (b) show the average TPR and PPV, and (c) shows the coefficient error $E_2$. The proposed Ident-WV (red) improves robustness compared to WeakIdent (blue) when the noise level is higher than 20%, showing a higher TPR and PPV, and a smaller coefficient error $E_2$ .

## 5.2 One-dimensional PDEs with initial conditions at various frequencies

We test Ident-WV on various one-dimensional (1D) differential equations listed in Table 1, where the initial condition varies with frequencies low $\omega = 2$, medium $\omega = 4$ and high $\omega = 8$ respectively. We perform 20 experiments of PDE identification with independent noise using Ident-WV, WeakIdent, WeakSINDy and compare their average performance using TPR, PPV for the support error, and $E_2$ for the coefficient error. The full results for all equations are presented in Figure 12, 13 and 14 in Appendix A. In Figure 5 and 6, we present the representative results. Figure 5 shows the results for the transport equation with diffusion and the Burgers equation when the initial condition has a low frequency: $\omega = 2$. In Figure 6, in the top row, we show the Burgers' equation where the initial condition has a medium frequency of $\omega = 4$, and in the bottom row, we present the Transport equation with diffusion where the initial condition has a high frequency $\omega = 8$. In general, when the initial condition has a low frequency, Ident-WV consistently stabilizes WeakIdent, and outperforms WeakSINDy. When the initial condition has medium and high frequencies, Ident-WV also stabilizes WeakIdent, and outperforms WeakSINDy in most cases.

| Equation | Initial Condition |
|---|---|
| **The Heat equation** $$u_t = 0.1592 u_{xx} \quad (30)$$ | $u(x,0) = \sin^2(\omega \pi x)$ |
| **The Transport equation** $$u_t = -u_x \quad (31)$$ | $u(x,0) = \sin^2(\omega \pi x / 0.7)$ if $x \in [0, 0.7]$ |
| **The Transport equation with diffusion** $$u_t = -10 u_x + u_{xx} \quad (32)$$ | $u(x,0) = \sin^3(\omega \pi (x-1))$ for $x \in [1, 2]$ |
| **The Burgers' equation** $$u_t = -u u_x \quad (33)$$ | $u(x,0) = 100 \sin(\omega \pi x)$ |
| **The Burgers' equation with diffusion** $$u_t = -u u_x + 0.2 u_{xx} \quad (34)$$ | $u(x,0) = 100 \sin(\omega \pi (x - 0.5))$ for $x \in [0.5, 1.5]$ |

Table 1: Various 1D equations and initial conditions with frequency $\omega$. The frequencies $\omega = 2, 4,$ and 8, are considered in experiments presented in Figure 5, 6, 12, 13 and 14.

## 5.3 1D equation with higher order Derivatives and 2D equations

We next present more experiments on 1D equations with higher order derivatives and 2D equations as listed in Table 2. These equations are more difficult to identify than the ones in Table 1, and we test the initial conditions in Table 2 without varying the frequencies.

Figure 7 shows the identification results for the Kuramoto–Sivashinsky (KS) equation (35). We evaluated 20 independent trials as the noise-to-signal ratio $\sigma_{\text{NSR}}$ varies from 0 to 150%. Figure 7 compares Ident-WV, WeakIdent, and WeakSINDy in terms of TPR, PPV and coefficient error. Across the entire noise range, Ident-WV and WeakIdent exhibit comparable performance in terms of TPR and PPV, indicating that both methods are effective at identifying the correct support of the governing equation. In contrast, WeakSINDy displays a distinct behavior: it achieves a relatively high TPR, even in high-noise regimes, implying that it tends to identify most of the true terms. However, this comes at the cost of a lower PPV, indicating a tendency to include spurious terms (false positives) in the model.

Figure 8 shows the identification results for the nonlinear Schrödinger (NLS) equation (36), based on 20 independent trials for each value of the noise-to-signal ratio $\sigma_{\text{NSR}} \in [0, 40\%]$. The NLS equation consists of two variables $u$ and $v$, and we use the collection of dynamic indicators $R[u]$, $R[u^2]$, $R[(u^2)_x]$, $R[(u^2)_{xx}]$, $R[(u^2)_t]$ to weight the equations about $u_t$, and the collection of dynamic indicators $R[v]$, $R[v^2]$, $R[(v^2)_x]$, $R[(v^2)_{xx}]$, $R[(v^2)_t]$ to weight the equations about $v_t$. The performance of three methods, the proposed Ident-WV, WeakIdent, and WeakSINDy, is evaluated in terms of support recovery (via TPR and PPV) and coefficient error ($E_2$).

In Figure 8, both Ident-WV and WeakIdent demonstrate high accuracy in support recovery across all noise levels. They maintain near-perfect True Positive Rate (TPR) and Positive Predictive Value (PPV) up to $\sigma_{\text{NSR}} \approx 20\%$, and only exhibit a gradual degradation as noise increases. The $E_2$ coefficient error remains low and comparable for both methods, indicating that they are equally effective in recovering the true model structure and parameter values under moderate noise. In contrast, WeakSINDy yields a relatively lower TPR when $\sigma_{\text{NSR}} \geq 20\%$, and suffers from a lower PPV, indicating frequent inclusion of false positive terms. Its $E_2$ coefficient error also grows more rapidly as noise level increases and exhibits higher variance, reflecting sensitivity in coefficient estimation when $\sigma_{\text{NSR}} \geq 20\%$. In this example, Ident-WV and WeakIdent are both robust in identifying the governing equation of the NLS equation under noise, with Ident-WV showing slightly better stability with smaller $E_2$ coefficient error at higher noise levels.

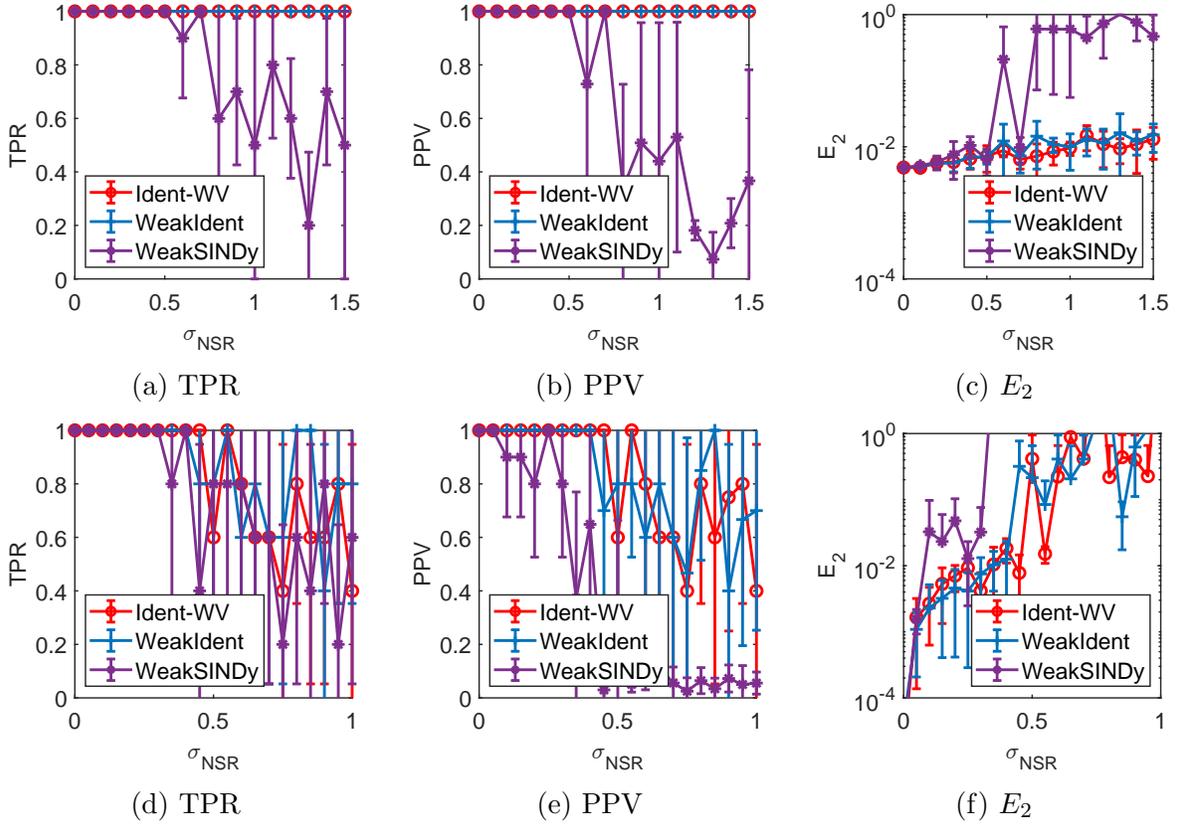

Figure 5: Identification examples of the 1D PDEs in Table 1, (the full results are presented in Figure 12) when the initial condition has a low frequency: $\omega = 2$. The top row is the Transport equation with diffusion (32), and the second row is the Burgers' equation (33). The $E_2$ coefficient error as shown a function of the noise-to-signal ratio $\sigma_{\text{NSR}}$. Note that in TPR and PPV graphs, Ident-WV (red) are closer to 1 and the $E_2$ error in the last column is smaller compared to WeakIdent (blue) and WeakSINDy (purple). When the initial condition has a low frequency, Ident-WV consistently stabilizes WeakIdent, and outperforms WeakSINDy.

Figure 9 shows the identification results for the PM equation (37). The PM equation is in 2D. For the dynamics-guided weighted weak form, we use the collection of dynamics indicators from the reference features that contain the spatial derivatives of the dependent variable $(u, u^2)$ up to second order, and its first-order time derivative. These weight matrices are $R[u]$, $R[u^2]$, $R[(u^2)_x]$, $R[(u^2)_y]$, $R[(u^2)_{xx}]$, $R[(u^2)_{xy}]$, $R[(u^2)_{yy}]$, $R[(u^2)_t]$. In Figure 9, (a) - (c) display the average TPR, PPV and $E_2$ coefficient error when NSR $\sigma_{\text{NSR}}$ increases from 0 to 7%. In this example, Ident-WV achieves higher TPR and PPV than WeakIdent and WeakSINDy at all levels, and it gives rise to a smaller coefficient error.

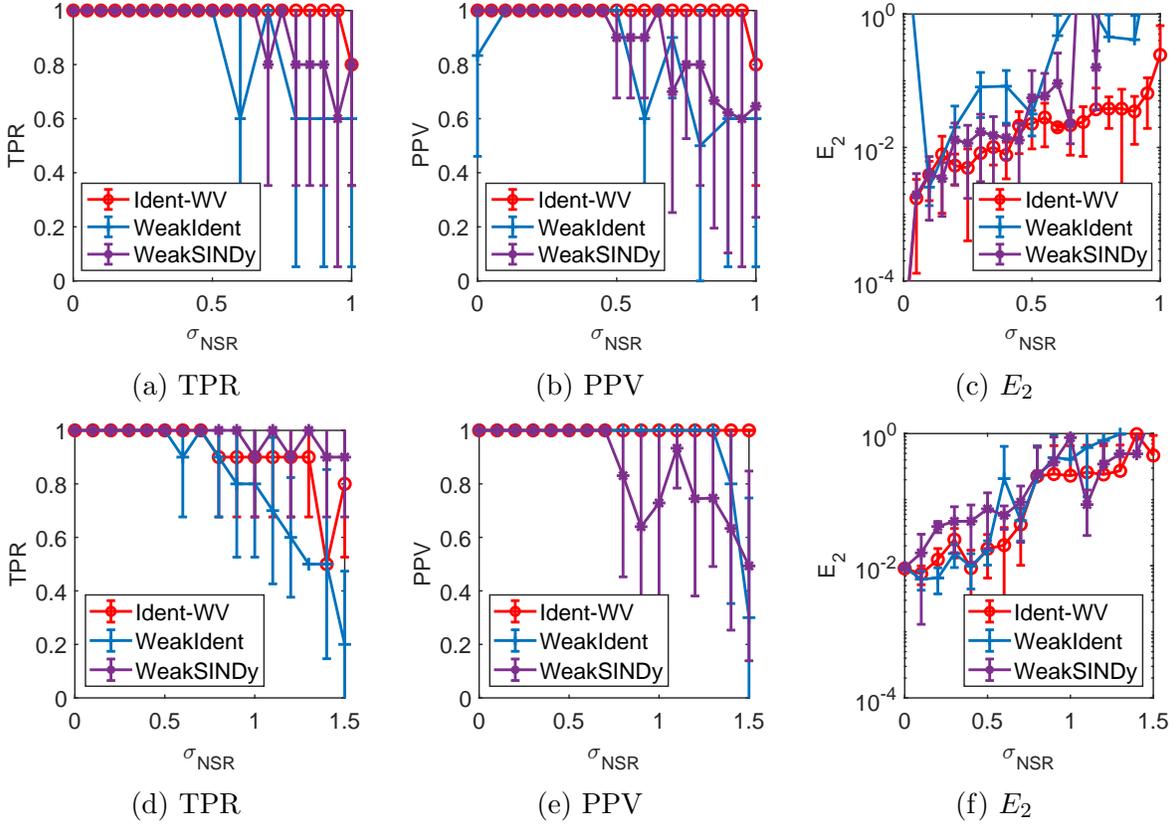

Figure 6: Identification examples of the 1D PDEs in Table 1 (the full results are presented in Figure 13 and 14). The top row shows Burgers' equation (33) with a medium frequency $\omega = 4$, and the bottom row shows the transport equation with diffusion (32) with a high frequency $\omega = 8$. In TPR and PPV graphs, Ident-WV (red) are closer to 1 and the $E_2$ error in the forth column is smaller compared to WeakIdent (blue) and WeakSINDy (purple). Ident-WV stabilizes WeakIdent, and outperforms WeakSINDy in most cases.

### 5.4 Comparison with Ensemble-SINDy

In this section, we compare our proposed Ident-WV on the Burgers' equation (33) with the Ensemble-SINDy (E-SINDy) method [8]. E-SINDy incorporates WeakSINDy [17] and two ensemble techniques: library bagging and row bagging. In library bagging, one randomly subsamples a collection of features from the library, and then identifies a differential equation from each subsampled feature library. The bagging step aggregates the equations identified from all subsampled feature libraries. In row bagging, one randomly subsamples some rows from the feature matrix, and then identifies a differential equation from each subsample. The bagging step aggregates the equations identified from all subsampled rows.

To systematically assess the impact of each ensemble technique, we test E-SINDy under four configurations: 1) E-SINDy-B: using only row bagging with 100 ensembles; 2) E-SINDy-LB: using only library bagging; 3) E-SINDy-DB-5: applying both bagging steps with 5 row ensembles; 4) E-SINDy-DB-100: applying both bagging steps with 100 row ensembles. Our experiments on these four ensemble techniques are given in Appendix B, in Figure 15 for Burgers' equation, and Figure 16 for KS equation.

| Equation | Initial Condition |
|---|---|
| **Kuramoto-Sivashinsky (KS)** $$u_t = -uu_x - u_{xx} - u_{xxxx} \quad (35)$$ | $u(x,0) = \cos(x/16)(1 + \sin(x/16))$ |
| **Nonlinear Schrodinger (NLS)** $$\begin{cases} u_t = 0.1u_{xx} + 0.1u_{yy} + u + v^3 - uv^2 + u^2v - u^3 \\ v_t = 0.1v_{yy} + 0.1v_{xx} + v - v^3 - uv^2 - u^2v - u^3 \end{cases} \quad (36)$$ | $u(x,0) = 2\operatorname{sech}(x)$ $v(x,0) = 0$ |
| **Anisotropic Porous Medium (PM)** $$u_t = 0.3u_{yy} - 0.8u_{xy} + u_{xx} \quad (37)$$ | $u(x,y,0) = \max(-0.26786x^2 - 0.71429xy - 0.89286y^2 + 0.4611, 0)$ |

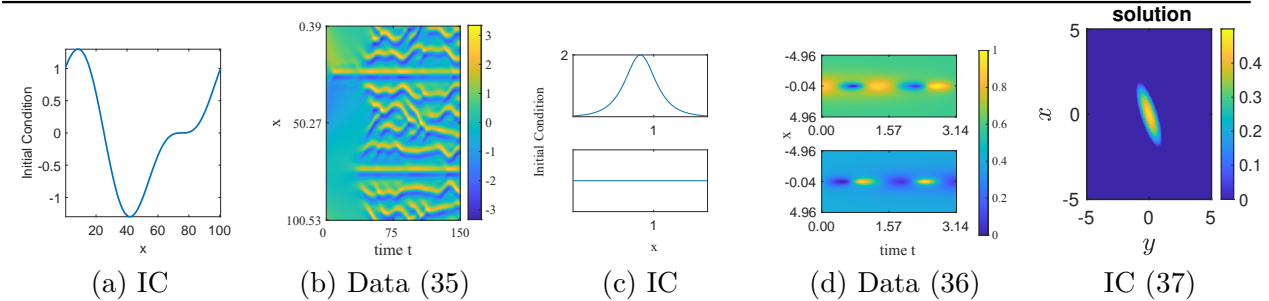

(a) IC  (b) Data (35)  (c) IC  (d) Data (36)  IC (37)

Table 2: Higher order 1D and 2D equations and their initial conditions considered in the experiments presented in Figure 7, 8 and 9. (a) and (b) are initial condition and the given data generated from the KS equation in (35), (c) and (d) are that of the Nonlinear Schrodinger equation in (36), and (e) is the initial condition for the Anisotropic Porous Medium equation in (37).

We observe that the best-performing technique in E-SINDy employs double bagging, and most performance gain is attributed to library bagging. Thus, we present comparison our method with E-SINDy-DB-5 in Figure 10 and 11. For the Burgers' equation in Figure 10, E-SINDy-DB-5 and IDENT-WV are comparable in support recovery, which give rise to higher TPR and PPV than WeakIdent and WeakSINDy. The coefficient error of Ident-WV is slightly lower than that of E-SINDy-DB-5. For the KS equation in Figure 11, E-SINDy-DB-5 yields slightly higher TPR than Ident-WV, with the sacrifice of a lower PPV.

## 6  Conclusion

In this paper, we used the dynamics-guided weighted weak form to develop Ident-WV to identify differential equations from a single trajectory of noisy data. We proposed the dynamics-guided weight matrix to highlight the high dynamic region of a reference feature and a voting strategy to stabilize the identification results from several reference features. Comprehensive numerical experiments are presented for Ident-WV, in comparison with WeakIdent, WeakSINDy and Ensemble-SINDy. Our experiments show that Ident-WV improves WeakIdent, and demonstrates superior robustness to noise across all equations tested in this paper. Our results underscore the effectiveness of Ident-WV in accurately identifying the underlying equation from noisy observations. Its ability to combine the robustness of weak formulations with enhanced stability from support and coefficient voting makes it particularly suitable for real-world scientific data, where noise is unavoidable and accurate model recovery is critical. Overall, this study highlights the promise of weak form-based approaches for PDE discovery and motivates the further development of noise-tolerant, data-efficient methods for interpretable scientific machine learning.

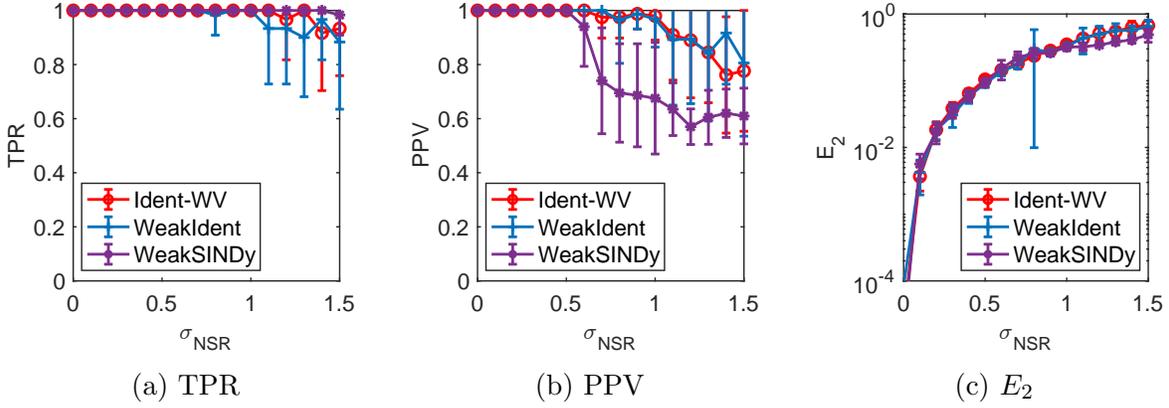

Figure 7: Identification results for the KS equation (35) in Table 2. (a) shows the initial condition and (b) the given data of the KS equation. (c), (d) and (e) , respectively, show the average TPR, PPV and $E_2$ coefficient error when the noise-to-signal ratio $\sigma_{\text{NSR}}$ increases from 0 to 150%. The standard deviation in these 20 experiments is represented by the vertical bar.

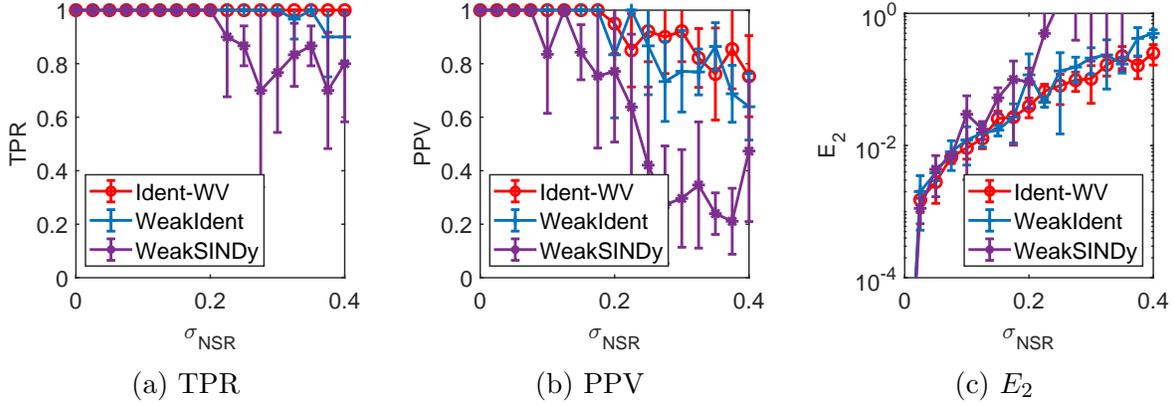

Figure 8: Identification results for the NLS equation (36) in Table 2. In (a) (b) (c), the average TPR, PPV and $E_2$ coefficient error are shown when $\sigma_{\text{NSR}}$ increases from 0 to 40%. The standard deviation in these 20 experiments is represented by the vertical bar.

## Acknowledgment

Wenjing Liao is supported by the National Science Foundation under the NSF DMS 2145167 and the U.S. Department of Energy under the DOE SC0024348. Haomin Zhou is partially supported by the NSF DMS 2307465.

## References

[1] Ellen Baake, Michael Baake, HG Bock, and KM Briggs. Fitting ordinary differential equations to chaotic data. *Physical Review A*, 45(8):5524, 1992.

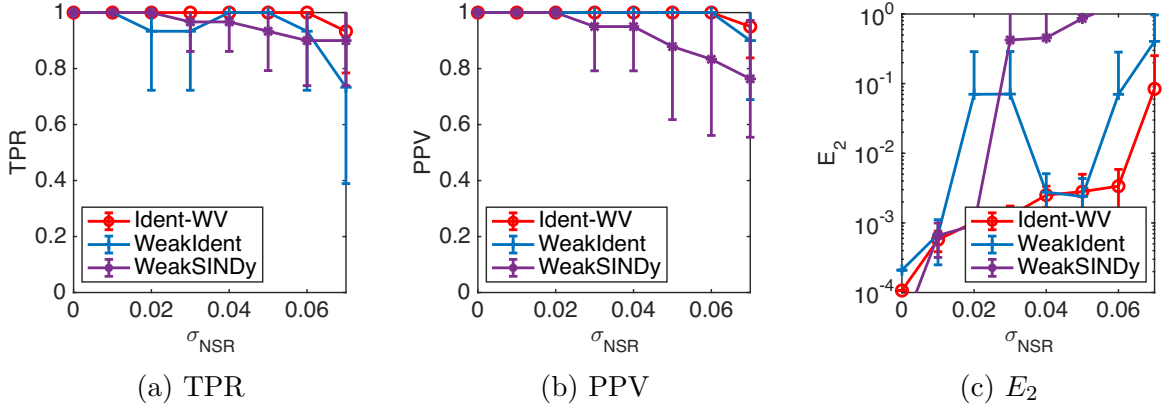

(a) TPR      (b) PPV      (c) $E_2$

Figure 9: Identification results for the PM equation (37) in Table 2. (a) - (c) display the average TPR, PPV and $E_2$ coefficient error when NSR $\sigma_{\text{NSR}}$ increases from 0 to 7%. The standard deviation in these 20 experiments is represented by the vertical bar.

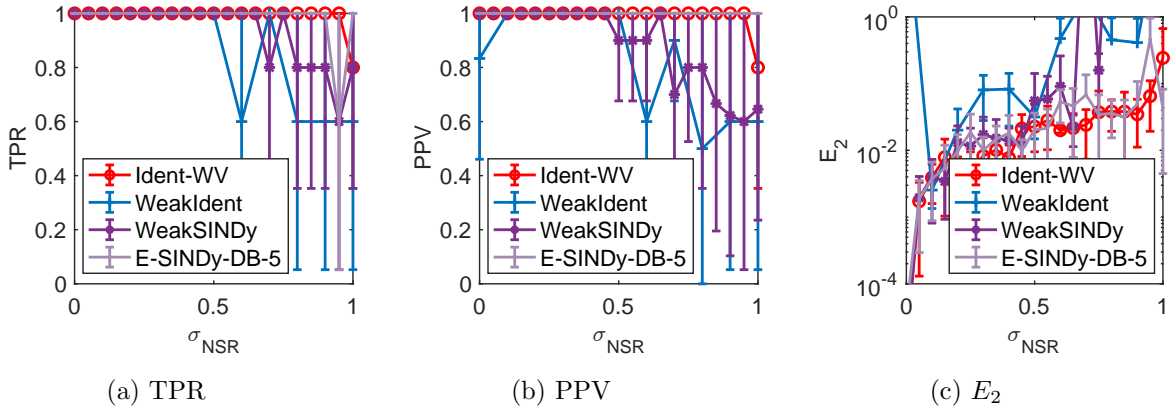

(a) TPR      (b) PPV      (c) $E_2$

Figure 10: [Comparison with an ensemble method] Identification results on the Burgers' equation (33) with the initial condition $u(x,0) = 100\sin(4\pi x)$ by Ident-WV, WeakIdentm WeakSINDy and E-SINDy-DB-5. Our proposed Ident-WV and E-SINDy-DB-5 give rise to the highest TPR and PPV. The coefficient error of Ident-WV is slightly lower than that of E-SINDy-DB-5.

[2] Markus Bär, Rainer Hegger, and Holger Kantz. Fitting partial differential equations to space-time dynamics. *Physical Review E*, 59(1):337, 1999.

[3] Hans Georg Bock. Recent advances in parameter identification techniques for ode. *Numerical treatment of inverse problems in differential and integral equations*, pages 95–121, 1983.

[4] Josh Bongard and Hod Lipson. Automated reverse engineering of nonlinear dynamical systems. *Proceedings of the National Academy of Sciences*, 104(24):9943–9948, 2007.

[5] Steven L Brunton, Joshua L Proctor, and J Nathan Kutz. Discovering governing equations from data by sparse identification of nonlinear dynamical systems. *Proceedings of the National Academy of Sciences*, 113(15):3932–3937, 2016.

[6] Emmanuel J Candes and Terence Tao. Decoding by linear programming. *IEEE transactions on information theory*, 51(12):4203–4215, 2005.

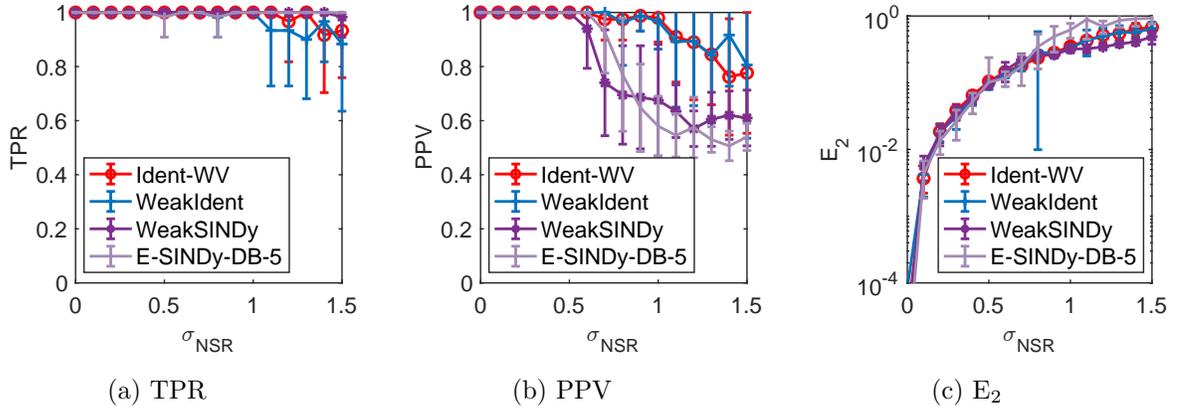

(a) TPR  (b) PPV  (c) $E_2$

Figure 11: [Comparison with an ensemble method] Identification results on the KS equation (35) by Ident-WV, WeakIdent, WeakSINDy and E-SINDy-DB-5. The initial condition is the same as provided in Table 2. In this example, E-SINDy-DB-5 yields slightly higher TPR than Ident-WV, with the sacrifice of a lower PPV.

[7] David L Donoho, Xiaoming Huo, et al. Uncertainty principles and ideal atomic decomposition. *IEEE transactions on information theory*, 47(7):2845–2862, 2001.

[8] Urban Fasel, J Nathan Kutz, Bingni W Brunton, and Steven L Brunton. Ensemble-sindy: Robust sparse model discovery in the low-data, high-noise limit, with active learning and control. *Proceedings of the Royal Society A*, 478(2260):20210904, 2022.

[9] Daniel R Gurevich, Patrick AK Reinbold, and Roman O Grigoriev. Robust and optimal sparse regression for nonlinear pde models. *Chaos: An Interdisciplinary Journal of Nonlinear Science*, 29(10):103113, 2019.

[10] Roy Y He, Haixia Liu, and Hao Liu. Group projected subspace pursuit for block sparse signal reconstruction: Convergence analysis and applications. *Applied and Computational Harmonic Analysis*, 75:101726, 2025.

[11] Yuchen He, Sung Ha Kang, Wenjing Liao, Hao Liu, and Yingjie Liu. Numerical identification of nonlocal potential in aggregation. *Communications in Computational Physics*, 2022.

[12] Yuchen He, Sung-Ha Kang, Wenjing Liao, Hao Liu, and Yingjie Liu. Robust identification of differential equations by numerical techniques from a single set of noisy observation. *SIAM Journal on Scientific Computing*, 44(3):A1145–A1175, 2022.

[13] Yuchen He, Sung Ha Kang, Wenjing Liao, Hao Liu, and Yingjie Liu. Group projected subspace pursuit for identification of variable coefficient differential equations (gp-ident). *Journal of Computational Physics*, 494:112526, 2023.

[14] Yuchen He, Hongkai Zhao, and Yimin Zhong. How much can one learn a partial differential equation from its solution? *Foundations of Computational Mathematics*, 24(5):1595–1641, 2024.

[15] Sung Ha Kang, Wenjing Liao, and Yingjie Liu. Ident: Identifying differential equations with numerical time evolution. *Journal of Scientific Computing*, 87(1):1–27, 2021.

[16] Lennart Ljung. *System identification*. Springer, 1998.

[17] Daniel A Messenger and David M Bortz. Weak sindy for partial differential equations. *Journal of Computational Physics*, page 110525, 2021.

[18] Daniel A Messenger and David M Bortz. Weak sindy: Galerkin-based data-driven model selection. *Multiscale Modeling & Simulation*, 19(3):1474–1497, 2021.

[19] Thomas G Müller and Jens Timmer. Parameter identification techniques for partial differential equations. *International Journal of Bifurcation and Chaos*, 14(06):2053–2060, 2004.

[20] Patrick AK Reinbold, Daniel R Gurevich, and Roman O Grigoriev. Using noisy or incomplete data to discover models of spatiotemporal dynamics. *Physical Review E*, 101(1):010203, 2020.

[21] Hayden Schaeffer. Learning partial differential equations via data discovery and sparse optimization. *Proceedings of the Royal Society A: Mathematical, Physical and Engineering Sciences*, 473(2197):20160446, 2017.

[22] Michael Schmidt and Hod Lipson. Distilling free-form natural laws from experimental data. *Science*, 324(5923):81–85, 2009.

[23] Cheng Tang, Roy Y. He, and Hao Liu. Wg-ident: Weak group identification of pdes with varying coefficients, 2025.

[24] Mengyi Tang, Wenjing Liao, Rachel Kuske, and Sung Ha Kang. Weakident: Weak formulation for identifying differential equation using narrow-fit and trimming. *Journal of Computational Physics*, 483:112069, 2023.

[25] Mengyi Tang, Hao Liu, Wenjing Liao, and Sung Ha Kang. Fourier features for identifying differential equations (fourierident). *arXiv preprint arXiv:2311.16608*, 2023.

[26] Kailiang Wu and Dongbin Xiu. Numerical aspects for approximating governing equations using data. *Journal of Computational Physics*, 384:200–221, 2019.

## A  Full results in Section 5.2

In this section, we demonstrate the full results about 1D equations in Table 1 in Section 5.2 when the initial conditions have various frequencies.

Figure 12 illustrates the identification results for the 1D PDEs in Table 1 when the initial condition has a low frequency: $\omega = 2$. The first, second and third rows, respectively, show the TPR, PPV and the $E_2$ coefficient error as a function of the noise-to-signal ratio $\sigma_{\text{NSR}}$. When the initial condition has a low frequency, Ident-WV consistently stabilizes WeakIdent, and outperforms WeakSINDy.

Figure 13 shows the identification of the 1D PDEs in Table 1 when the initial condition has a medium frequency: $\omega = 4$. The first, second and third rows respectively show the TPR, PPV and the $E_2$ coefficient error as a function of the noise-to-signal ratio $\sigma_{\text{NSR}}$. In comparison with Figure 12, the identification results are improved for all methods When the initial condition has a medium frequency. In most cases, Ident-WV stabilizes WeakIdent, and outperforms WeakSINDy.

Figure 14 shows the identification of the 1D PDEs in Table 1 when the initial condition has a high frequency: $\omega = 8$. Still Ident-WV stabilizes WeakIdent, and outperforms WeakSINDy in most cases.

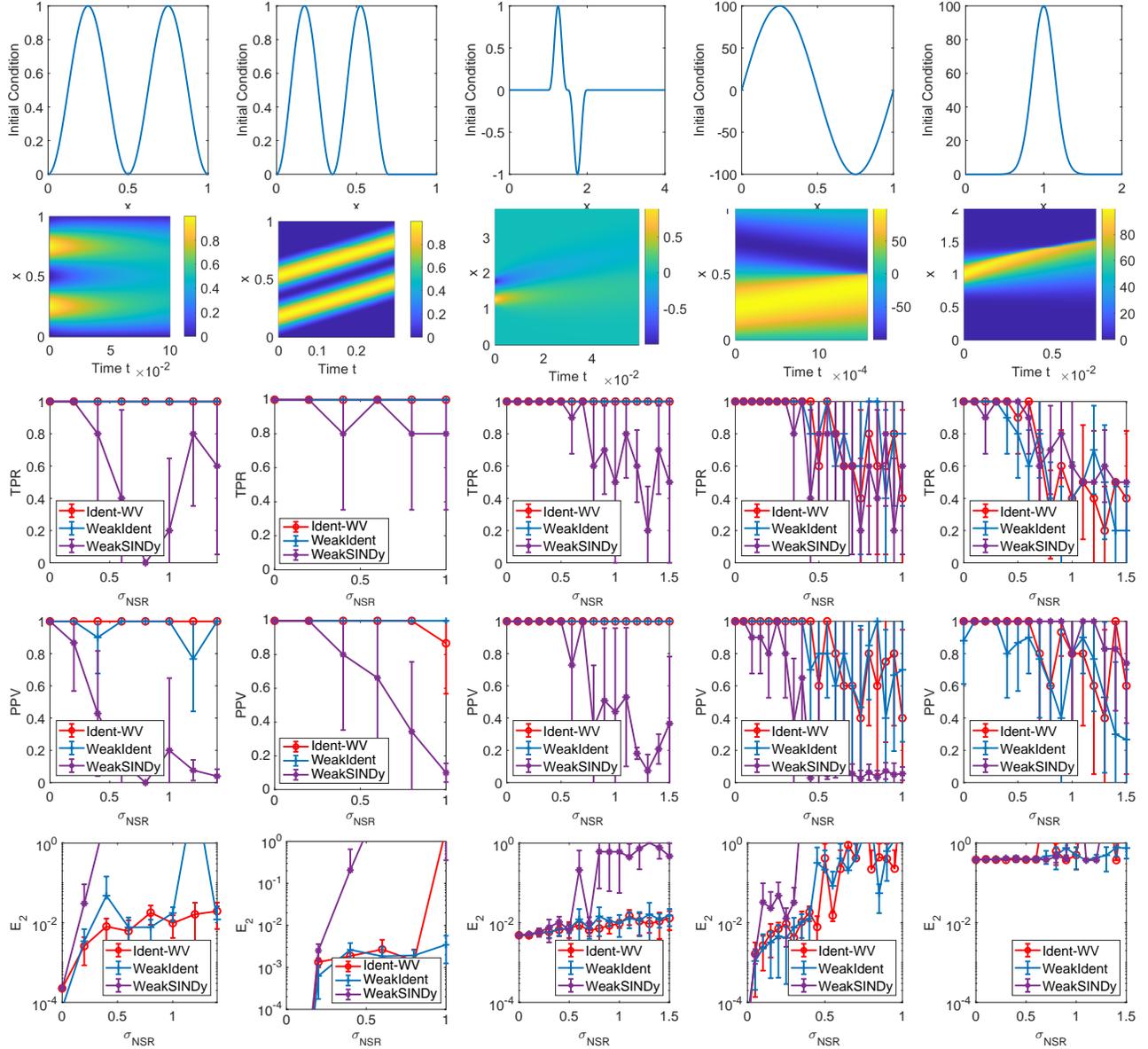

(a) Heat equation  (b) Transport Equation  (c) Transport with Diffusion  (d) Burgers' Equation  (e) Burgers' Equation with diffusion

Figure 12: Identification of the 1D PDEs in Table 1 when the initial condition has a low frequency: $\omega = 2$. Each column contains the result for one equation in Table 1. The first, second and third rows respectively show the TPR, PPV and the $E_2$ coefficient error as a function of the noise-to-signal ratio $\sigma_{\text{NSR}}$. When the initial condition has a low frequency, Ident-WV consistently stabilizes WeakIdent, and outperforms WeakSINDy.

## B  A full comparison with Ensemble-SINDy

We test Ensemble-SINDy (E-SINDy) [8] under four configurations: (a) E-SINDy-B: using only row bagging with 100 ensembles; (b) E-SINDy-LB: using only library bagging; (c) E-SINDy-DB-5: applying both bagging steps with 5 row ensembles; (d) E-SINDy-DB-100: applying both bagging

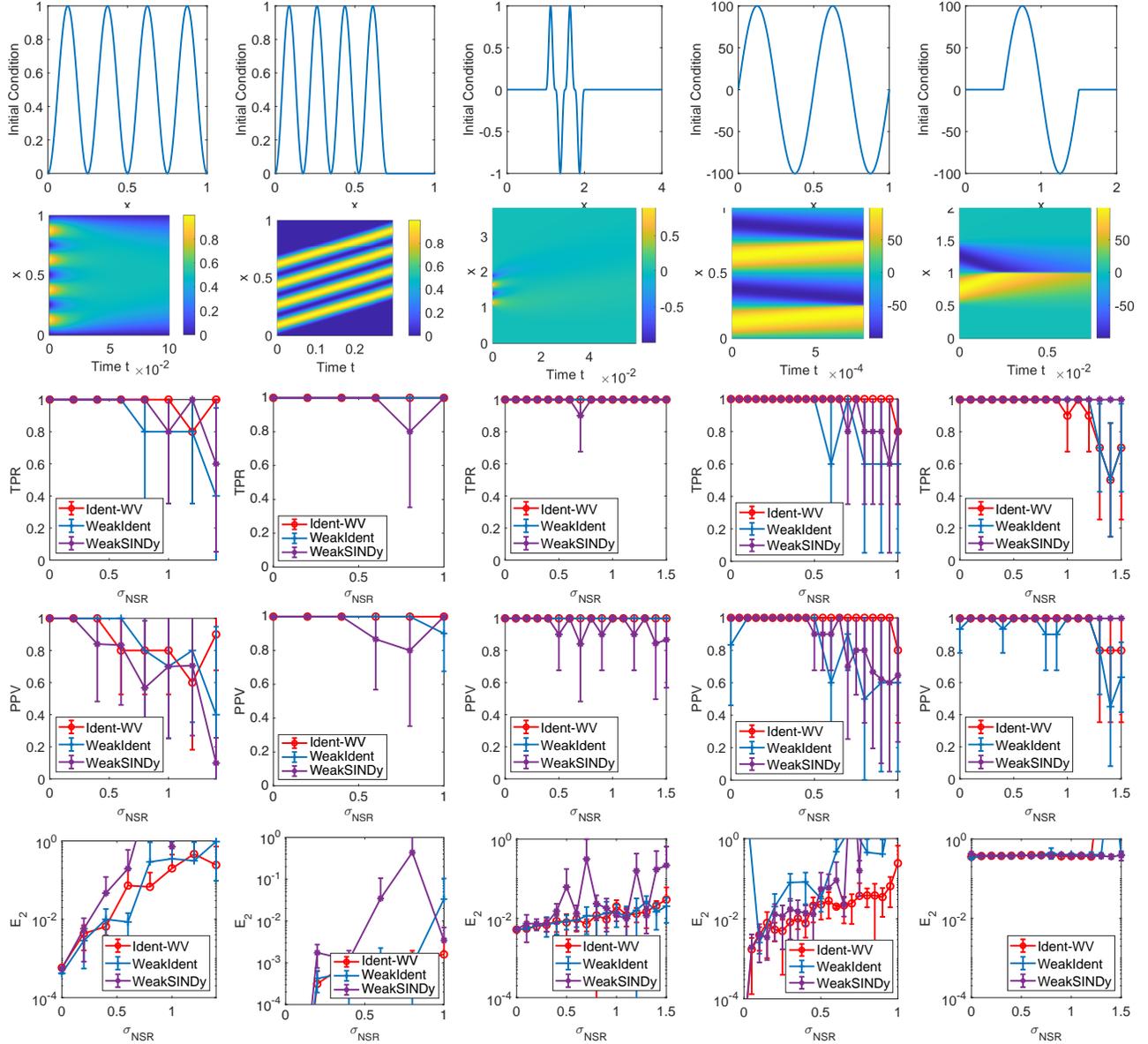

(a) Heat equation  (b) Transport Equation  (c) Transport with Diffusion  (d) Burgers' Equation  (e) Burgers' Equation with diffusion

Figure 13: Identification of the 1D PDEs in Table 1 when the initial condition has a medium frequency: $\omega = 4$. Each column contains the result for one equation in Table 1. The first, second and third rows respectively show the TPR, PPV and the $E_2$ coefficient error as a function of the noise-to-signal ratio $\sigma_{\text{NSR}}$. In comparison with Figure 12, the identification results are improved for all methods When the initial condition has a medium frequency. Ident-WV stabilizes WeakIdent, and outperforms WeakSINDy in most cases.

steps with 100 row ensembles. Figure 15 shows a comparison for the Burgers' equation (33), and Figure 16 shows that for the KS equation (35). These experiments demonstrate that the major improvement in E-SINDy is achieved through library bagging. For the Burgers' equation with analytic solution, E-SINDy has comparable performance to our proposed Ident-WV. However, for

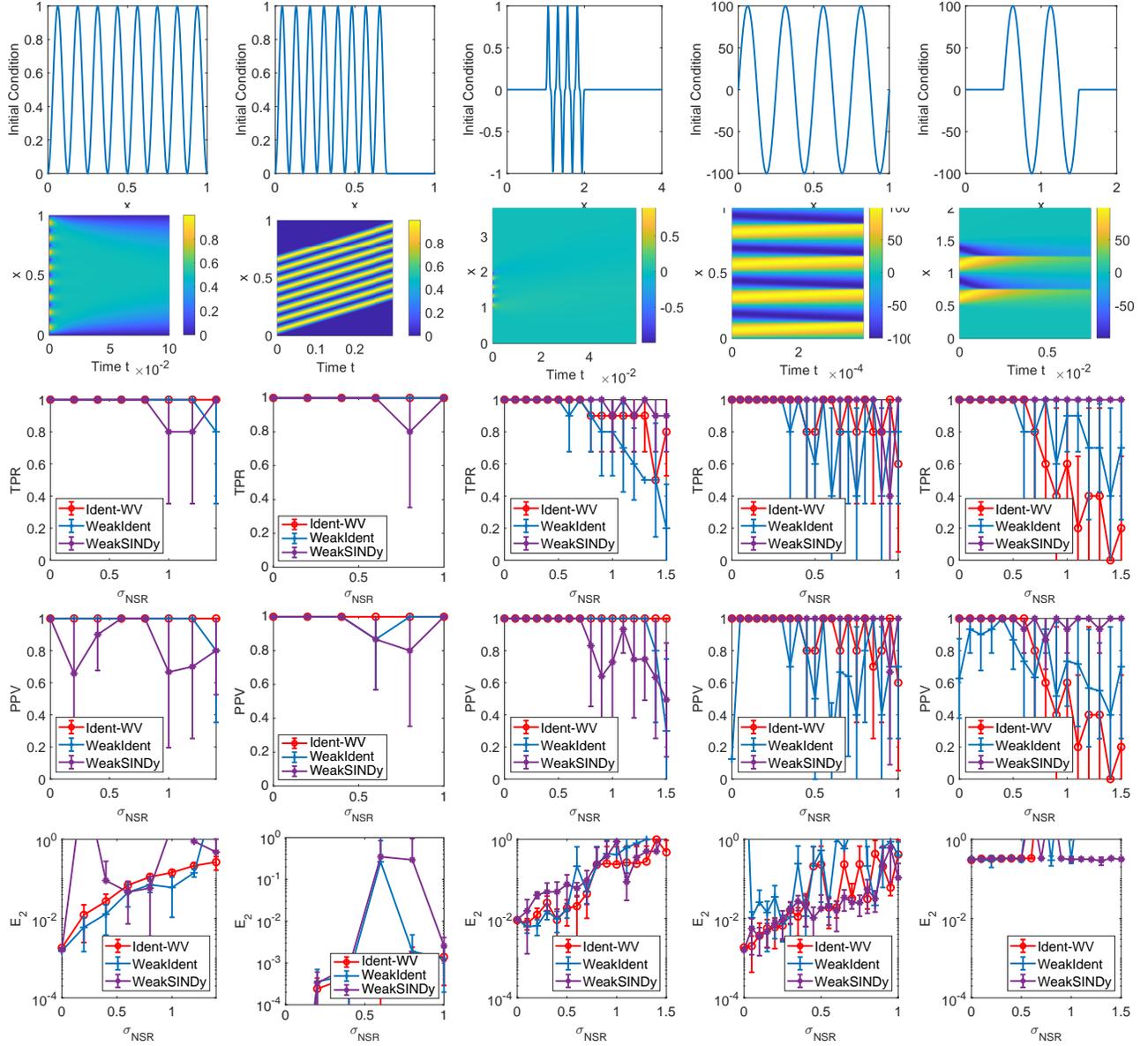

(a) Heat equation  (b) Transport Equation  (c) Transport with Diffusion  (d) Burgers' Equation  (e) Burgers' Equation with diffusion

Figure 14: Identification of the 1D PDEs in Table 1 when the initial condition has a high frequency: $\omega = 8$. Each column contains the result for one equation in Table 1. The first, second and third rows respectively show the TPR, PPV and the $E_2$ coefficient error as a function of the noise-to-signal ratio $\sigma_{\text{NSR}}$. Still Ident-WV stabilizes WeakIdent, and outperforms WeakSINDy in most cases.

the KS equation, under high noise levels, the PPV of E-SINDy remains lower than that of Ident-WV.

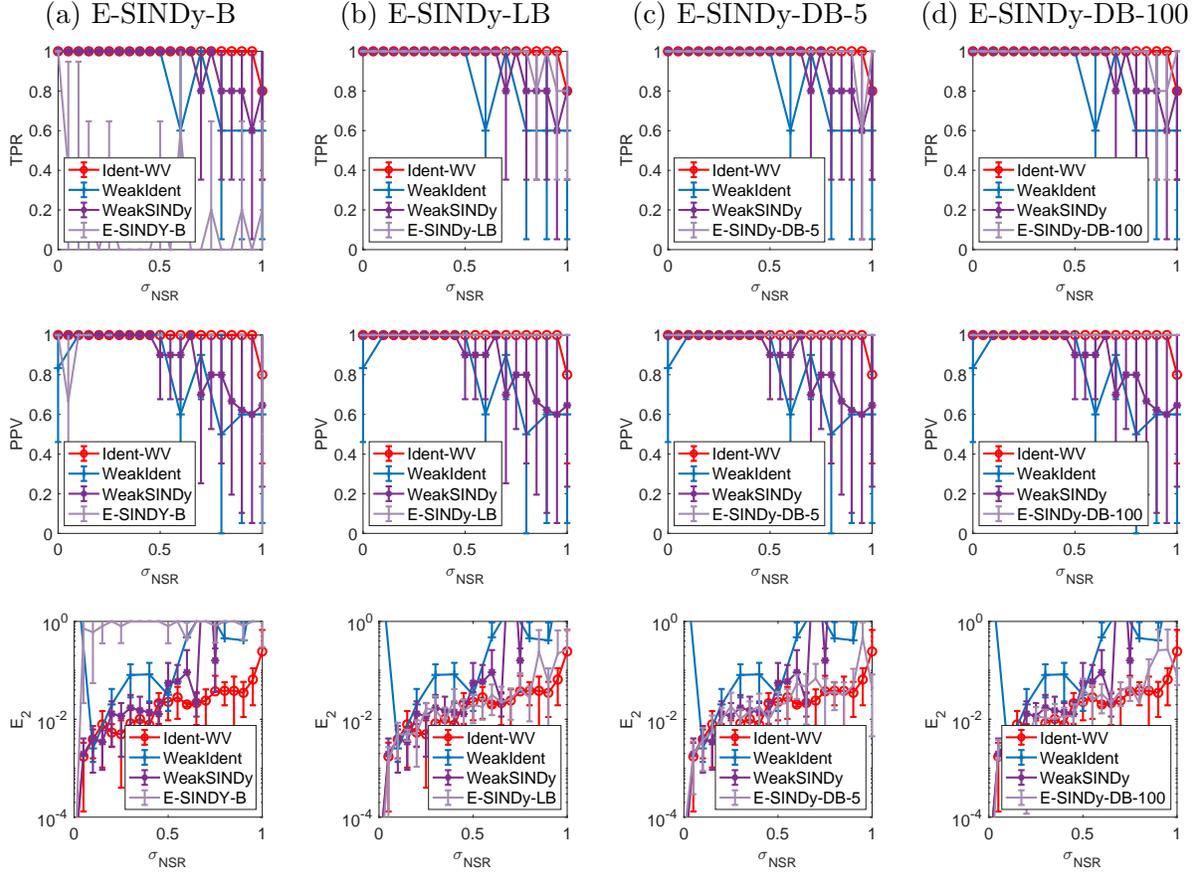

Figure 15: Comparison of our proposed method Iden-WV and Ensemble-SINDy with multiple configurations (a) E-SINDy-B, (b) E-SINDy-LB, (c) E-SINDy-DB-5 and (d) E-SINDy-DB-100. They are compared on the Burgers' equation (33) with the initial condition $u(x,0) = 100\sin(4\pi x)$. The first, second and third rows respectively show the TPR, PPV and the $E_2$ coefficient error as a function of the noise-to-signal ratio $\sigma_{\text{NSR}}$.

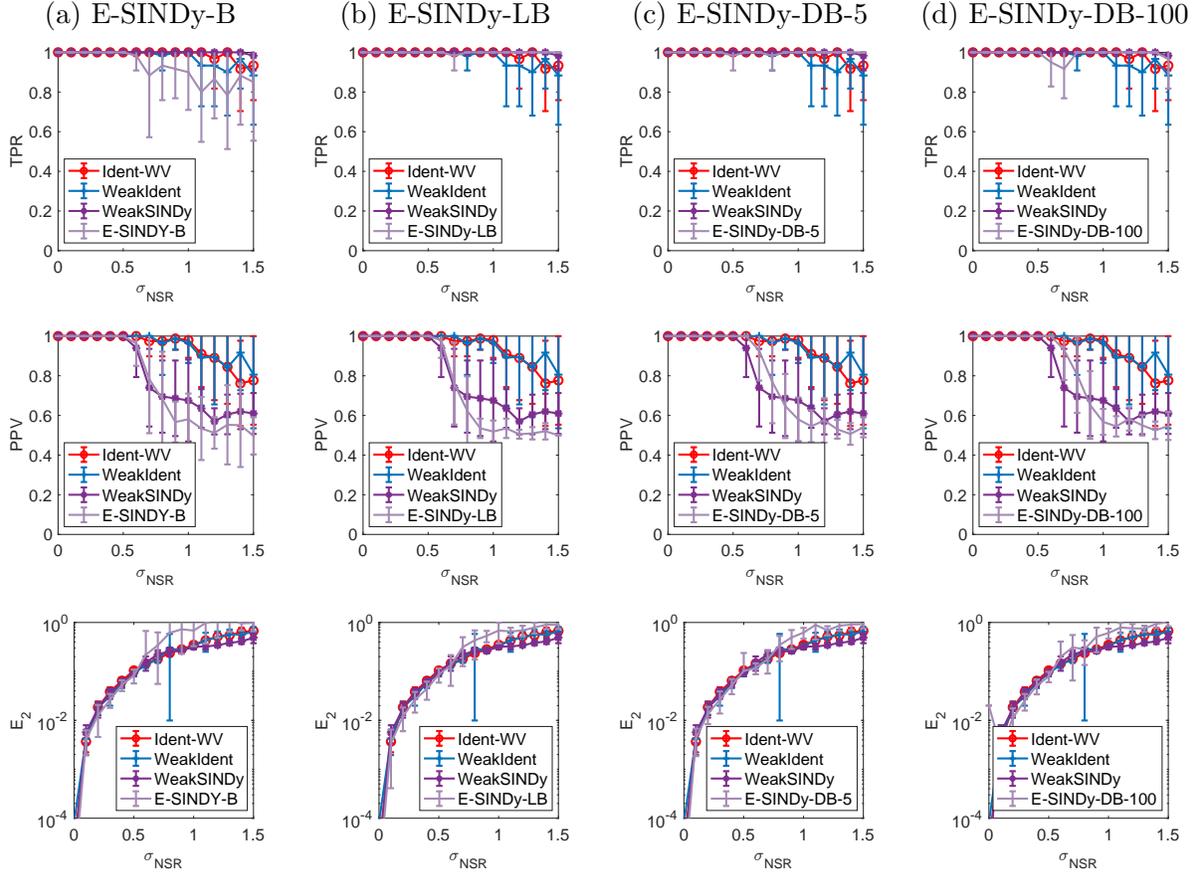

Figure 16: Comparison of our proposed method Iden-WV and Ensemble-SINDy with multiple configurations: (a) E-SINDy-B, (b) E-SINDy-LB, (c) E-SINDy-DB-5 and (d) E-SINDy-DB-100. They are compared on the KS equation (35) with the initial condition presented in Table 2. . Each column contains the result for one configuration mentioned in Section 5.4. The first, second and third rows respectively show the TPR, PPV and the $E_2$ coefficient error as a function of the noise-to-signal ratio $\sigma_{\text{NSR}}$.